\documentclass[11pt]{article}

\usepackage{amsmath,amssymb,amsthm}
\usepackage{enumerate,color,bbm,ifthen,graphicx,overpic}
\usepackage{mathrsfs}
\usepackage{a4wide}
\usepackage[latin1]{inputenc}
\usepackage{psfrag}
\usepackage{graphics}
\usepackage{subfigure}

\newtheorem{theo}{Theorem}[section]
\newtheorem{lemma}[theo]{Lemma}

\newtheorem{prop}[theo]{Proposition}

\theoremstyle{remark}
\newtheorem{rem}[theo]{Remark}

\numberwithin{equation}{section}

\def\N{\mathbb{N}}
\def\t{\theta}

\def\eqsp{\;}
\def\Xset{\mathsf{X}} 
\def\Xsigma{\mathcal{X}} 

\def\PP{\mathbb{P}} 

\def\barX{\overline{X}}
\newcommand{\eqdef}{\ensuremath{\stackrel{\mathrm{def}}{=}}}

\newcommand{\un}{{\mathbbm{1}}}

\DeclareMathAlphabet{\mathpzc}{OT1}{pzc}{m}{it}
\newcommand{\dps}{\displaystyle } 
\newcommand{\rme}{\mathrm{e}} 
\newcommand{\EE}{\mathbb{E}}
\newcommand{\G}{\mathpzc{Geo}}
\newcommand{\eps}{\varepsilon}

\newcounter{hypoconbis}
\newcounter{saveconbis}
\newcommand\debutA{\begin{list} {\textbf{A\arabic{hypoconbis}}}{\usecounter{hypoconbis}}\setcounter{hypoconbis}{\value{saveconbis}}}
\newcommand\finA{\end{list}\setcounter{saveconbis}{\value{hypoconbis}}}
\newenvironment{adem}[1][]%
   {\ \\ {\bf Proof#1. }}%
   {\hfill\mbox{\rule{2 true mm}{3 true mm}}\vskip 2 ex\noindent}

\begin{document}

\title{Efficiency of the Wang-Landau algorithm: a simple test case}
\author{ G. Fort$^1$, B. Jourdain$^2$, E. Kuhn$^3$, T Leli\`evre$^{2,4}$ and G. Stoltz$^{2,4}$ \\
\footnotesize{1: LTCI, CNRS \& Telecom ParisTech, 46 rue Barrault, 75634 Paris Cedex 13, France} \\
\footnotesize{2: Universit\'e Paris Est, CERMICS, Ecole des Ponts ParisTech, 
6 \& 8 Av. Pascal, 77455 Marne-la-Vall\'ee Cedex 2, France}\\
\footnotesize{3: INRA Unité MIA, Domaine de Vilvert, 78352 Jouy-en-Josas Cedex, France}\\
\footnotesize{4: Project-team MICMAC, INRIA Rocquencourt, Domaine de Voluceau, B.P. 105, 78153 Le Chesnay Cedex, France}  
}

\maketitle

\begin{abstract}
We analyze the efficiency of the Wang-Landau algorithm to sample
a multimodal distribution on a prototypical simple test case. We show that
the exit time from a metastable state is much smaller for the Wang
Landau dynamics than for the original standard Metropolis-Hastings
algorithm, in some asymptotic regime. Our results are confirmed by
numerical experiments on a more realistic test case.
\end{abstract}

\section{Introduction}

The Wang-Landau algorithm was originally proposed in the physics literature to
efficiently sample the density of states of Ising-type
systems~\cite{wang-landau-01,WL01PRL}. It belongs to the class of {\em free energy biasing
techniques}~\cite{lelievre-rousset-stoltz-07-b} which have been introduced in
computational statistical physics to efficiently sample thermodynamic
ensembles and to compute free energy differences. From a computational statistical
point of view, it can be seen as some adaptive importance sampling strategy
combined with a Metropolis algorithm~\cite{MRRTT53,Hastings70}: the target probability distribution is updated
at each iteration of the algorithm in order to have a sampling of the
configuration space as uniform as possible along a given direction. There are numerous physical and biochemical
works using this technique to overcome sampling problems such as the ones
encountered in the computation of macroscopic properties around critical points
and phase transitions, or for the sampling of folding mechanisms for proteins. The original
paper~\cite{WL01PRL} is cited more than one thousand times, according to Web of Knowledge.
The success of the technique motivated its use and study in the statistics
literature, see~\cite{Liang05,LLC07,atchade-liu-10,JR11,BJDD11,FJKLS1} for instance
for previous mathematical and numerical studies.

There are in fact several variations of the original Wang-Landau
algorithm, see the discussion in~\cite{FJKLS1}. We study here the
Wang-Landau algorithm with a deterministic adaption sequence (see
Section~\ref{sec:linearized_WL} for a precise definition of the
algorithm). The aim of this article is to discuss from a mathematical
viewpoint the efficiency of the Wang-Landau procedure. The real
practical interest of adaptive importance sampling techniques is indeed their improved convergence properties, compared to standard sampling techniques. Although this improvement is obvious to practitioners, it is mathematically more difficult to formalize.

This paper is a companion paper to~\cite{FJKLS1} where a convergence
result is proven, without any efficiency analysis. Actually, to our
knowledge, the previous mathematical studies on the Wang-Landau algorithm
solely focused on the convergence of the algorithm, not on its
efficiency. Such insight into improved convergence rates has been
obtained for other adaptive importance sampling methods, in particular
for Adaptive Biasing Force techniques,
see~\cite{LRS08,lelievre-minoukadeh-09}. These analysis have been
performed on the nonlinear Fokker Planck equation obtained in the
limit of infinitely many interacting replicas. To the best of our
knowledge, there is currently no efficiency analysis of adaptive importance sampling
techniques based on a {\em single} trajectory interacting with its own
past. The aim of this work is to gain some insight on the efficiency of
the Wang-Landau algorithm, which is an example of such a technique. More precisely,
we
show here through the analytical study of a toy model and a confirmation by
numerical results in a more complicated case, that the Wang-Landau algorithm
indeed allows to efficiently escape from metastable states.

\medskip 

The paper is organized as follows. We describe in Section~\ref{sec:description_algorithm} the algorithm that we consider.  We next turn to a discussion on the efficiency of the method in Section~\ref{sec:efficiency}. On a very simple example, we mathematically quantify the improvement on the convergence properties given by the Wang-Landau dynamics, compared to a standard Metropolis Hastings procedure.  Our results are confirmed by numerical experiments on a more realistic two-dimensional test case presented in Section \ref{sec:exit_numerical}. The proofs of our results are gathered in Section~\ref{sec:proof:cvg}. Section \ref{sec:succtimes} is devoted to some refinement of the comparison between the standard Metropolis Hastings procedure and the Wang-Landau algorithm.

\section{Description of the Wang-Landau algorithm}
\label{sec:description_algorithm}

\subsection{Notation and preliminaries}

Let us consider a normalized target probability density
$\pi$ defined on a Polish space $\Xset$, endowed with a reference
measure~$\lambda$ defined on the Borel $\sigma$-algebra $\Xsigma$. As for classical Metropolis-Hastings procedure, the practical implementation of the
algorithm only requires to specify $\pi$ up to a multiplicative
constant.
 In statistical physics, the set $\Xset$ is typically composed of all
admissible configurations of the system while $\pi$ is a Gibbs measure with
density $\pi(x) = Z_\beta^{-1} \exp(-\beta U(x))$, $U$ being the potential
energy function and $\beta$ the inverse temperature.  In condensed matter
physics for instance, actual simulations are performed on systems composed of
$N$ particles in dimension 2 or 3, living in a cubic box with periodic
boundary conditions. In this case, $\Xset = (L\mathbb{T})^{2N}$ or
$\Xset = (L\mathbb{T})^{3N}$, where $L$ is the length of the sides of
the box and
$\mathbb{T} = \mathbb{R} / \mathbb{Z}$ is the one-dimensional torus.

We now consider a partition $\Xset_1, \ldots,
\Xset_d$ of $\Xset$ in $d \geq 2$ elements, and define, for any $i \in \{1, \ldots, d \}$,
\begin{equation}
  \label{eq:def:thetastar}
 \t_\star(i) \eqdef \int_{\Xset_i} \pi(x) \lambda(dx) \eqsp.
\end{equation}
In the following, $\Xset_i$ will be called the $i$-th {\em stratum}.  Each
weight $\t_\star(i)$ is assumed to be positive and gives the relative
likelihood of the stratum $\Xset_i \subset \Xset$.  In practice, the
partitioning could be obtained by considering some smooth function $\xi \, : \,
\Xset \to [a,b]$ (called a reaction coordinate in the physics literature) and
defining, for $i=1,\dots,d-1$,
\begin{equation}
  \label{eq:def_Xi}
\Xset_i = \xi^{-1}\Big([\alpha_{i-1},\alpha_{i})\Big) \eqsp,
\end{equation}
and $\Xset_d =  \xi^{-1}\left([\alpha_{d-1},\alpha_{d}]\right)$, with 
$a = \alpha_0 < \alpha_1 < \dots <\alpha_d = b$ (possibly, $a=-\infty$
and/or $b=+\infty$). 

Let us emphasize here that the choice of an appropriate function $\xi$ is a difficult but central issue. It is mostly based on intuition at the time being: practitioners 
identify some slowly evolving degrees of freedom
responsible for the metastable behavior of the system, and build a
function $\xi$ and then a partition using these slow degrees of
freedom. Here, metastability refers to the fact that trajectories
generated by the reference (non-adaptive) dynamics, which is ergodic with
respect to the target probability measure $\pi$ (for example a Metropolis Hastings algorithm with target $\pi$), remain trapped for a long time 
in some region of $\Xset$, and only occasionally hop to
another region, where they also remain trapped. There are ways to 
quantify the relevance of the choice of the reaction coordinate, see for instance
the discussion in~\cite{CLS12}. There are also ways to adaptively choose the levels $(\alpha_i)_{0\leq i\leq d}$, see~\cite{BJDD11}.

The above discussion motivates the fact that the weights $\t_\star(i)$
typically span several orders of magnitude, some sets $\Xset_i$ having very
large weights, and other ones being very unlikely under~$\pi$. Besides,
trajectories bridging two very likely states typically need to go through unlikely
regions.  To efficiently explore the configuration space, and sample numerous
configurations in all the strata $\Xset_i$, it is therefore a natural idea to resort
to importance sampling strategies and to appropriately reweight each subset~$\Xset_i$. A possible way to do so is the following.  Let $\Theta$
be the subset of (non-degenerate) probability measures on $\{1, \ldots, d \}$ given by
\[
\Theta = \left\{
\t = (\t(1), \ldots, \t(d)) \ \left| \ 0 < \t(i) <1 \ \text{for all $i \in
  \{1, \ldots, d \}$ and} \ \sum_{i=1}^d \t(i) =1 \right. \right\} \eqsp.
\]
For any $\t \in \Theta$, we denote by $\pi_\t$ the probability density on $(\Xset,
\Xsigma)$ (endowed with the reference measure~$\lambda$) defined as
\begin{equation}
  \label{eq:def:pitheta}
 \pi_\t(x) = \left(\sum_{i=1}^d \frac{\t_\star(i)}{\t(i)}\right)^{-1} \ 
\sum_{i=1}^d \frac{\pi(x)}{\t(i)} \ \un_{\Xset_i}(x) \eqsp. 
\end{equation}
This measure is such that the weight of the set $\Xset_i$ under $\pi_\t$ is
proportional to $\t_\star(i)/\t(i)$. In particular, all the strata $\Xset_i$
have the same weight under~$\pi_{\t_\star}$. Unfortunately, the vector $\t_\star$ is
unknown and sampling under $\pi_{\theta_\star}$ is typically
unfeasible. 

The Wang-Landau algorithm allows precisely to overcome these difficulties: at each iteration of the algorithm, a
weight vector $\t_n = (\t_n(1), \ldots, \t_n(d))$ is updated based on the past
behavior of the algorithm and a new point is drawn from a Markov kernel $P_{\t_n}$
with invariant density $\pi_{\t_n}$. The update of $\{\t_n, n \geq 0\}$
is chosen in such a way to penalize the already visited strata.
The intuition for the convergence of
this algorithm is that if $\{\t_n, n\geq 0 \}$ converges to $\t_\infty$ then the
draws are asymptotically distributed according to the density
$\pi_{\t_\infty}$ and it can be checked from the updating rule (see
Equation~\eqref{eq:NL_update} below) that necessarily $\t_\infty=\t_\star$.



\subsection{The Wang-Landau algorithm with  deterministic
  adaption}
\label{sec:linearized_WL}

 We now
describe the algorithm we study in this article.  Let $\{\gamma_n, n \geq 1\}$
be a $[0,1)$-valued deterministic sequence. For any $\t \in \Theta$, denote by
$P_\t$ a Markov transition kernel onto $(\Xset, \Xsigma)$ with unique
stationary distribution $\pi_\t(x)\lambda(dx)$; for example, $P_\t$ is
one step of a Metropolis-Hastings
algorithm~\cite{MRRTT53,Hastings70} with target probability measure $\pi_\t(x)\lambda(dx)$.

Consider an initial value $X_0 \in \Xset$ and an initial set of weights $\t_0
\in \Theta$ (typically, in absence of any prior information, $\t_0(i) = 1/d$).
Define the process $\{(X_n, \t_n), n\geq 0 \}$ as follows: given the current
value $(X_n, \t_n)$,
\begin{enumerate}[\quad (1)]
\item Draw $X_{n+1}$ under the conditional
  distribution $P_{\t_n}(X_n, \cdot)$;
\item The weights are then updated as
\begin{equation}
\label{eq:NL_update}
\t_{n+1}(i) = \t_n(i) \frac{1 + \gamma_{n+1}
  \un_{\Xset_i}(X_{n+1})}{\dps 1 + \gamma_{n+1}\theta_n(I(X_{n+1}))}
\text{  for all } i \in \{1, \ldots, d\} \eqsp.
\end{equation}
Here, $I: \Xset \to \{1, \ldots,
d\}$ defined by
\begin{equation}
  \label{eq:definition:fonctionI}
  \forall x \in \Xset, \, I(x)=i \text{ if and only if } x \in \Xset_i
\end{equation}
associates to a point $x$ the index $I(x)$ of the stratum where $x$ lies.
\end{enumerate}
As explained above, the idea of the updating
strategy~\eqref{eq:NL_update} is that the weights of the
visited strata are increased, in order to penalize already visited states. Note that the update~\eqref{eq:NL_update} is such that the sum of the weights remains equal to~1.

Let us recall the result of convergence proved in~\cite{FJKLS1}. Three
assumptions are required: on the equilibrium
measure (see~A\ref{hyp:targetpi}), on the transition kernels $\{P_\t, \t
\in\Theta \}$ (see~A\ref{hyp:kernel}) and on the step-size sequence
$\{\gamma_n, n \geq 1\}$ (see~A\ref{hyp:stepsize}). It is assumed that
\debutA
\item \label{hyp:targetpi} The probability density $\pi$ with respect to the
  measure $\lambda$ is such that $0 < \inf_{\Xset}\pi \leq \sup_\Xset \pi <
  \infty$.  In addition, $\inf_{1 \leq i \leq d } \lambda(\Xset_i)> 0$. \finA
Notice that Assumption~A\ref{hyp:targetpi} implies that
$\inf_{1 \leq i \leq d } \t_\star(i)> 0$ where
  $\t_\star$ is given by~(\ref{eq:def:thetastar}). 
\debutA
  \item \label{hyp:kernel} For any $\t \in \Theta$, $P_\t$ is a
    Metropolis-Hastings transition kernel with invariant distribution $\pi_\t \ 
    d\lambda$, where $\pi_\t$ is given by (\ref{eq:def:pitheta}), and with symmetric
    proposal kernel $q(x,y) \lambda(dy)$ satisfying $\inf_{\Xset^2} q >
    0$.
\finA 

\debutA
\item \label{hyp:stepsize} The sequence $\{\gamma_n, n \geq 1\}$ is a
  $[0,1)$-valued deterministic sequence such that 
  \begin{enumerate}[a)]
  \item \label{hyp:pasdecroiss} $\{\gamma_n, n \geq 1 \}$ is a (ultimately)
    non-increasing sequence and $\lim_n \gamma_n = 0$;
  \item \label{hyp:stepsize:item2} $\sum_n \gamma_n  =\infty$; 
  \item \label{hyp:stepsize:item3} $\sum_n \gamma_n^2 < \infty$.
   \end{enumerate}
   \finA
A typical choice for the step-size sequence $\{\gamma_n,\, n \geq 1\}$ is
$\gamma_n=\gamma_\star n^{-\alpha} \eqsp$,
with $1/2 < \alpha \leq 1$.

Under assumptions A\ref{hyp:targetpi}-A\ref{hyp:kernel}-A\ref{hyp:stepsize}, it is shown in~\cite{FJKLS1} that the
algorithm converges: 
$$\PP\left( \lim_{n\to+\infty} \t_n = \t_\star \right) = 1.$$
More precisely, the proof is done for a slightly different update
than~\eqref{eq:NL_update}, namely the following linearized version:
\begin{equation}
\label{eq:update_weights_linearized}
\left\{ \begin{array}{ll}
    \t_{n+1}(i) =   \t_n(i) + \gamma_{n+1} \ \t_n(i)  \left(1 -  \t_n(i)\right),  & \\
    \t_{n+1}(k) = \t_n(k) - \gamma_{n+1} \ \t_n(k) \ \t_n(i), & \text{for} \ k \neq i.
\end{array} \right.
\end{equation}
The update~\eqref{eq:update_weights_linearized} is obtained
from~\eqref{eq:NL_update} in the limit of small $\gamma_n$. We believe
that the arguments used in~\cite{FJKLS1} can be adapted to prove the
convergence for the nonlinear update~\eqref{eq:NL_update}. By contrast, we would like to emphasize here that the distinction between the two
updating strategies~\eqref{eq:update_weights_linearized}
and~\eqref{eq:NL_update} {\em does matter} when
considering the flat histogram criterium for the update of the step-sizes, as proved in~\cite{JR11}.

However, this convergence result does not help to understand the success of the
Wang-Landau algorithm. This algorithm is actually known to be useful in
metastable situations, namely when the original Markov chain (with
transition kernel $P_{\t_0}$) remains trapped for very long times in some
regions (called the metastable states). Metastability is one of the
major bottleneck of standard Markov Chain Monte Carlo techniques,
since ergodic averages should be considered over very long times in
order to obtain accurate results. The aim of this article is to show
that in such a metastable situation, the Wang-Landau algorithm indeed
is an efficient sampling procedure. Our analysis will be based on
 estimates of exit times from
metastable states.

\section{Analytical results in a simple case}
\label{sec:efficiency}
We present in this section results on the improved convergence properties of
the Wang-Landau algorithm (when compared to non-adaptive samplers), by
theoretically analyzing the first exit times out of a
metastable state.  Indeed, adaptive biasing techniques such as the Wang-Landau
algorithm have been especially designed to be able to switch as fast as
possible from a metastable state to another in order to efficiently explore the
whole configuration space.

We show in this section that the Wang-Landau algorithm allows to
rapidly escape from a metastable state, namely from a large probability stratum
surrounded by small probability strata. More precisely,  we consider a
toy model composed of only three strata: two large probability
strata (the metastable states) separated by a low probability stratum (the transition
state). We are able to precisely quantify
the time the system needs to go from the first metastable state to the second
one, for adaptive and non-adaptive dynamics. We show in particular that the
exit time is dramatically reduced with the Wang-Landau dynamics compared to
the corresponding non-adaptive dynamics.

Using the notation of the previous section, we have only three strata
and three states, and thus
$\Xset=\{1,2,3\}$ and $\Xset_i=\{i\}$ for
$i=1,2,3$.  Jumps are only allowed between neighboring states, namely
from $1$ to $\{1,2\}$, from $2$ to $\{1,2,3\}$ and from $3$ to $\{2,3\}$. Though being very simple, we believe that this toy model is
prototypical of a metastable dynamics. We will check numerically in the next
section that our conclusions on this simple test case are indeed also
valid for more complicated and more realistic situations.

\subsection{Definition of the dynamics}\label{sec:def_dyn}

We assume that the first and third strata are visited with high probability, and that the second stratum is visited with
low probability. More precisely, we set
\begin{equation}\label{eq:theta_star}
\theta_\star(2) = \frac{\eps}{2+\eps} \eqsp, \qquad \theta_\star(1) =
\theta_\star(3) = \frac{1}{2+\eps} \eqsp,
\end{equation}
for a small positive parameter~$\eps\in (0,1)$, and consider the
limit $\eps \to 0$.  The target density $\pi$ on $\Xset$ is thus
defined as: $\pi(\{i\})=\theta_\star(i)$ for $i=1,2,3$ (the reference
measure $\lambda$ being the uniform measure on $\Xset=\{1,2,3\}$). The parameters $\theta_\star(i)$ depend on $\eps$, even
though we do not explicitly indicate this dependence to keep the notation
simple. In this specific setting, the biased probability
measure~\eqref{eq:def:pitheta} is
\[
\pi_\theta(i)=\left(\sum_{j=1}^3 \frac{\theta_\star(j)}{\theta(j)}
\right)^{-1} \frac{\theta_\star(i)}{\theta(i)} \text{ for } i \in \{1,2,3\}.
\]
Notice that $\pi_{\theta_\star}=(1/3,1/3,1/3)$ is the uniform measure on~$\Xset$.
 
The basic building block for the reference non-adaptive Markov chain 
$\{\barX_n, n\geq 0 \}$ is a symmetric proposal kernel allowing transitions to nearest-neighbor strata only:
\[
Q=\left[
\begin{array}{ccc}
\dps \frac{2}{3} & \dps \frac{1}{3} & 0 \\ [10pt]
\dps \frac{1}{3} & \dps \frac{1}{3} & \dps \frac{1}{3} \\ [10pt]
0     & \dps \frac{1}{3} & \dps \frac{2}{3}
\end{array}
\right] \eqsp.
\]
The corresponding non-adaptive Markov chain is built using a Metropolis
Hastings algorithm~\cite{MRRTT53,Hastings70}, with $Q$ as the proposal kernel, and $\pi$ as the
target distribution. To compute the kernel $\overline{P}$ of the
Metropolis algorithm, we evaluate its off-diagonal terms, and adjust
the diagonal in order for the rows to sum up to~1. For symmetric
proposals, the Metropolis procedure consists in proposing a new
configuration $\widetilde{X}_{n+1}$ from the previous state
$\overline{X}_{n}$ according to the proposal kernel~$Q$, and then to
accept this proposal with probability $1 \wedge \left(
\pi(\widetilde{X}_{n+1})/\pi(\overline{X}_{n}) \right)$, in which case
$\overline{X}_{n+1}=\widetilde{X}_{n+1}$; otherwise, $\overline{X}_{n+1} = \overline{X}_{n}$. For instance, the probability to go from~1 to~2 reads
\begin{equation}
\label{eq:P12}
\overline{P}_{12} = Q_{12} \left( 1 \wedge \frac{\pi(\{ 2 \})}{\pi(\{ 1 \})} \right) = \frac13  \left( 1 \wedge \frac{\theta_\star(2)}{\theta_\star(1)} \right).
\end{equation}
Since $\eps < 1$, the kernel $\overline{P}$ is given by
\begin{equation}
\label{eq:overline_P}
\overline{P}=\left[
\begin{array}{ccc}
\dps 1-\frac{\eps}{3} & \dps \frac{\eps}{3} & 0 \\[10pt]
\dps \frac{1}{3} & \dfrac{1}{3}  &  \dfrac{1}{3} \\[10pt]
0     &  \dfrac{\eps}{3} & 1-\dfrac{\eps}{3}  
\end{array}
\right] \eqsp.
\end{equation}
The non-adaptive dynamics
$\{\barX_n, n\geq 0 \}$ is metastable, in the sense that the time to go from
the stratum $1$ to the stratum $3$
\[
\overline{T}_{1 \to 3} = \min \Big\{n \, : \, \barX_n =3 \text{
  starting from } \barX_0=1 \Big\}
\]
is very large, and more precisely of order $6/\eps$ (see
Proposition~\ref{lem:non-adapt} below). This is due to the fact that, in order to go
from 1 to 3, the chain has to visit the very low probability
transition state~2. This is a prototypical metastable dynamics
reminiscent of what happens along molecular dynamics trajectories: due to
the very high dimensional configuration space, only local moves are
allowed (otherwise they would be mostly rejected) and thus, it is
difficult to go from a very likely region to another one since they
are usually
separated by low probability zones.

\medskip

For the associated adaptive Wang-Landau dynamics $\{(X_n,
\theta_n), n \geq 0\}$, the transition kernel $P_{\theta_n}$ to go from $X_n$ to $X_{n+1}$
is the Metropolis Hastings kernel corresponding to
 the proposal kernel $Q$ and the target distribution $\pi_{\t_n}$. The expression of $P_{\t}$ is obtained with computations similar to the ones leading to the expression~\eqref{eq:overline_P} of the transition kernel of the non-adaptive dynamics. In fact, it suffices to replace $\pi$ by $\pi_{\theta}$ in equalities such as~\eqref{eq:P12}. More precisely,
\begin{equation}\label{eq:Ptheta}
P_\theta=\left[
\begin{array}{ccc}
\dps 1-\frac{1}{3} \left(\eps\frac{\theta(1)}{\theta(2)} \wedge 1 \right) & \dps \frac{1}{3} \left(\eps\frac{\theta(1)}{\theta(2)} \wedge 1 \right) & 0 \\[10pt]
\dps \frac{1}{3} \left(\frac{1}{\eps}\frac{\theta(2)}{\theta(1)}
  \wedge 1 \right) & \dps 1 - \frac{1}{3}\left(\frac{1}{\eps}\frac{\theta(2)}{\theta(1)}
  \wedge 1  +  \frac{1}{\eps}\frac{\theta(2)}{\theta(3)} \wedge 1
\right)  &  \dps \frac{1}{3} \left(\frac{1}{\eps}\frac{\theta(2)}{\theta(3)} \wedge 1 \right) \\[10pt]
0     &  \dps \frac{1}{3}
\left(\eps\frac{\theta(3)}{\theta(2)} \wedge 1 \right) & \dps 1- \frac{1}{3}
\left(\eps\frac{\theta(3)}{\theta(2)} \wedge 1 \right)
\end{array}
\right] \eqsp.
\end{equation}
In addition, the step-size sequence in~\eqref{eq:NL_update} is
\begin{equation}\label{eq:gamma}
\gamma_n=\gamma_\star n^{-\alpha} \eqsp,
\end{equation}
for a positive constant $\gamma_\star$, and a parameter $\alpha \in
[1/2,1]$ (note that we allow here the value $1/2$, see Remark~\ref{rem:rmk_cv}).

We start from initially equiprobable strata 
$\theta_0(1)=\theta_0(2)=\theta_0(3)=1/3$, 
so that $\pi_{\theta_0} = \pi$. 
Notice that the non-adaptive dynamics is simply the Markov chain
with transition kernel $P_{(1/3,1/3,1/3)}$. It can be obtained from the adaptive
dynamics by setting $\gamma_\star = 0$, in which case $\theta_n =
\theta_0 = (1/3,1/3,1/3)$  for all $n \geq 0$. As above
for the non-adaptive dynamics, we define the time to go
from the stratum $1$ to the stratum $3$ for the Wang-Landau dynamics as
\[
{T}_{1 \to 3} = \min\Big\{n \, : \, X_n=3 \text{ starting from } X_0=1\Big\} \eqsp.
\]

The aim of this section is to show that, in some sense to be made precise,
$\overline{T}_{1\to 3}$ is much larger than $T_{1 \to 3}$ {\em i.e.}  the Wang-Landau
dynamics is much less metastable than the corresponding non-adaptive dynamics. This is related
to the fact that, when the stochastic process $\{X_n, n\geq 0\}$ remains
stuck in the stratum $1$, this stratum gets more and more penalized
($\theta_n(1)$
increases, see (\ref{eq:NL_update})), so that a transition to
the stratum~2 becomes more and more favorable.  From the stratum~2, a jump to
the stratum~3 is then very likely. This is the bottom line of the whole adaptive
procedure: penalizing the already visited strata in order to explore very
quickly new regions.

\subsection{Precise statement on the exit times}
We now provide a precise statement on how the exit times $\overline{T}_{1\to3}$ and
$T_{1 \to 3}$ scale when $\eps$ goes to zero.  For the non-adaptive dynamics, it holds
(see Section~\ref{sec:proof_non-adapt} for the proof):

\begin{prop}\label{lem:non-adapt}
  The time $\overline{T}_{1 \to 3}$ scales like~$6/\eps$, in the following sense:
\begin{align}\label{eq:espTbar13}
& \frac{\varepsilon}{6} \EE\left(\overline{T}_{1 \to 3}\right)= 1
+ \frac{\varepsilon}{2} \mathop{\longrightarrow}_{\eps
\to 0} 1 \eqsp, \\
\label{eq:probTbar13}
 & \forall c \geq 0, \qquad \lim_{\eps\to 0}\PP\left(
  \frac{\varepsilon}{6} \overline{T}_{1\to 3}> c\right) = \rme^{-c} \eqsp.
\end{align}
\end{prop}
Eq.~(\ref{eq:probTbar13}) states that when $\eps \to 0$, $\eps
\, \overline{T}_{1\to 3}$ converges in distribution to an exponential random variable
with parameter $1/6$.

\medskip

Let us now consider the Wang-Landau dynamics~\eqref{eq:WL}. The following
result holds (see Section~\ref{sec:proof_analytical_adaptive} for the proof).

\begin{prop}
\label{lem:analytical_adaptive} 
Let $\gamma_\star$ and $\alpha$ be
the two constants defining the sequence $\gamma_n$, as
given by~(\ref{eq:gamma}). Let us assume that $\alpha \in [1/2,1]$, with $\gamma_\star < 1$
if $\alpha=1/2$.  
\begin{itemize}
\item In the case $\alpha \in [1/2,1)$, the random variables 
$\left( | \ln \eps|^{-1/(1-\alpha)} \, T_{1 \to 3}\right)_{\eps >
  0}$ converge in probability to
$\left(\frac{1-\alpha}{\gamma_\star}\right)^{1/(1-\alpha)}$ when $\eps$
goes to $0$
\item In the case $\alpha =1$, for any function $h$ such that $\dps \lim_{\eps \to 0}
h(\eps)=+\infty$\eqsp, 
\begin{equation}\label{eq:probT:bis_lemme}
\lim_{\eps \to 0} \PP\left( \frac{1}{h(\varepsilon)}<\eps^{1/(1+\gamma_\star)} T_{1 \to 3}  <h(\varepsilon) \right) =1\eqsp.
\end{equation}
\end{itemize}
\end{prop}

In the case $\alpha=1$, one should think of
functions $h$ going very slowly to infinity, so that the above result
essentially means that
\begin{equation}
   \label{equivt}\mbox{ as }\varepsilon\to 0,
\qquad 
T_{1 \to 3}\mbox{ scales like }\begin{cases}\dps \left(\frac{1-\alpha}{\gamma_\star}\right)^{1/(1-\alpha)}|\ln \eps|^{1/(1-\alpha)}\mbox{ if $\alpha\in [1/2,1)$,}\\[10pt]
      \eps^{-1/(1+\gamma_\star)}\mbox{ if $\alpha=1$.}
   \end{cases}  
\end{equation}In any case, the Wang-Landau algorithm is such that $T_{1 \to 3}$ is much smaller than $\overline{T}_{1 \to 3}$ in the limit $\eps \to 0$ (namely in metastable situations). 

Notice that at time $T_{1\to 3}$, the Wang Landau algorithm cannot go back immediately to state $2$. It still has to get rid of part of the initial metastability : in particular $\tilde{\theta}_{T_{1\to 3}}(2)\geq\tilde{\theta}_{T_{1\to 3}}(3)$ since state $2$ has been visited at least once before $T_{1\to 3}$ and the sequence of step-sizes is decreasing. As a consequence, the entry $(3,2)$ of the matrix $P_{\theta_{T_{1\to 3}}}$ which gives the probability for the algorithm to go back to state $2$ at time $T_{1\to 3}+1$ is smaller than $\frac{\eps}{3}$. Section \ref{sec:succtimes} is dedicated to a formal analysis of the scaling in terms of $\eps$ of the successive durations between a visit by the algorithm of one of the extremal states $1$ and $3$ and the next visit of the other extremal state. Some hint at the total time necessary to get rid of the metastability is also given.

\begin{rem}
\label{rem:rmk_cv}
Two points should be mentioned about the convergence result from Proposition~\ref{lem:analytical_adaptive}:
\begin{itemize}
   \item The convergence in probability in the case
$\alpha \in [1/2,1)$ is equivalent to: for all $C_a$ and $C_b$ such
that $0<C_a < \left(\frac{1-\alpha}{\gamma_\star}\right)^{1/(1-\alpha)} \!\!\!\!\!< C_b$,
\[
\lim_{\eps \to 0} \PP\left( T_{1 \to 3} \in  \left(C_a|\ln \eps|^{1/(1-\alpha)}  , C_b  |\ln \eps|^{1/(1-\alpha)}\right) \right) =1.
\]
According to Proposition~\ref{lem:T13}, this limiting probability is
still one with a lower bound slightly larger than $C_a|\ln \eps|^{1/(1-\alpha)}$. 
\item Notice that we obtain results on first exit
times also for $\alpha =1/2$, which is an excluded value to obtain the
almost sure convergence of the Wang-Landau algorithm (see assumptions A3 above).
\end{itemize}
\end{rem}

\section{Numerical illustrations}
\label{sec:exit_numerical}

The aim of this section is to show that (most of) the results obtained
for the very simple three-state model of
Section~\ref{sec:efficiency} are still valid  for a less simple
example inspired by target measures used in computational statistical
physics. In these numerical experiments, we also investigate the
behavior of the algorithm for values of $\alpha$ in the interval $(0,1/2]$, which are excluded values to prove the theoretical convergence of stochastic approximation procedures in general, and in particular of the Wang-Landau algorithm (see assumptions~A3). The interest of choosing $\alpha \in (0,1/2]$ is that the Wang Landau algorithm escapes must faster from metastable states. It is therefore easier to numerically investigate very large values of $\beta$.

Our aim is to study the behavior of the exit times out of a
metastable state as the temperature in the system goes to zero. The
temperature will thus play a role similar to the role of~$\varepsilon$ in the
Section~\ref{sec:efficiency} (see formula~\eqref{eq:eps_beta}
below, where $\beta$ is the inverse temperature).

\subsection{Presentation of the model and of the dynamics}

We consider the system based on the two-dimensional potential suggested in~\cite{PSLS03}. The state space is $\Xset = [-R,R] \times \mathbb{R}$ (with $R > 0$), and we denote by $x=(x_1,x_2)$ a generic element of~$\Xset$. The reference measure $\lambda$ is the Lebesgue measure. The density of the target measure reads
\[
\pi(x) \propto \un_{[-R,R]}(x_1) \, \mathrm{e}^{-\beta U(x_1,x_2)},
\]
for some positive inverse temperature $\beta$, with
\begin{align}
U(x_1,x_2)
& = 3 \exp\left(-x_1^2 - \left(x_2-\frac13\right)^2\right)
- 3 \exp\left(-x_1^2 - \left(x_2-\frac53\right)^2\right) \label{eq:pot_U} \\
& \quad - 5 \exp\left(-(x_1-1)^2 - x_2^2\right)
- 5 \exp\left(-(x_1+1)^2 - x_2^2\right)
+ 0.2 x_1^4 + 0.2 \left(x_2-\frac13\right)^4. \nonumber
\end{align}
We introduce $d$ strata $\Xset_\ell = (a_\ell,a_{\ell+1}) \times
\mathbb{R}$, with $a_\ell = -R + 2(\ell-1) R/d$ and $\ell=1, \ldots,
d$. 

\begin{figure}\begin{center}
\includegraphics[width=7.4cm]{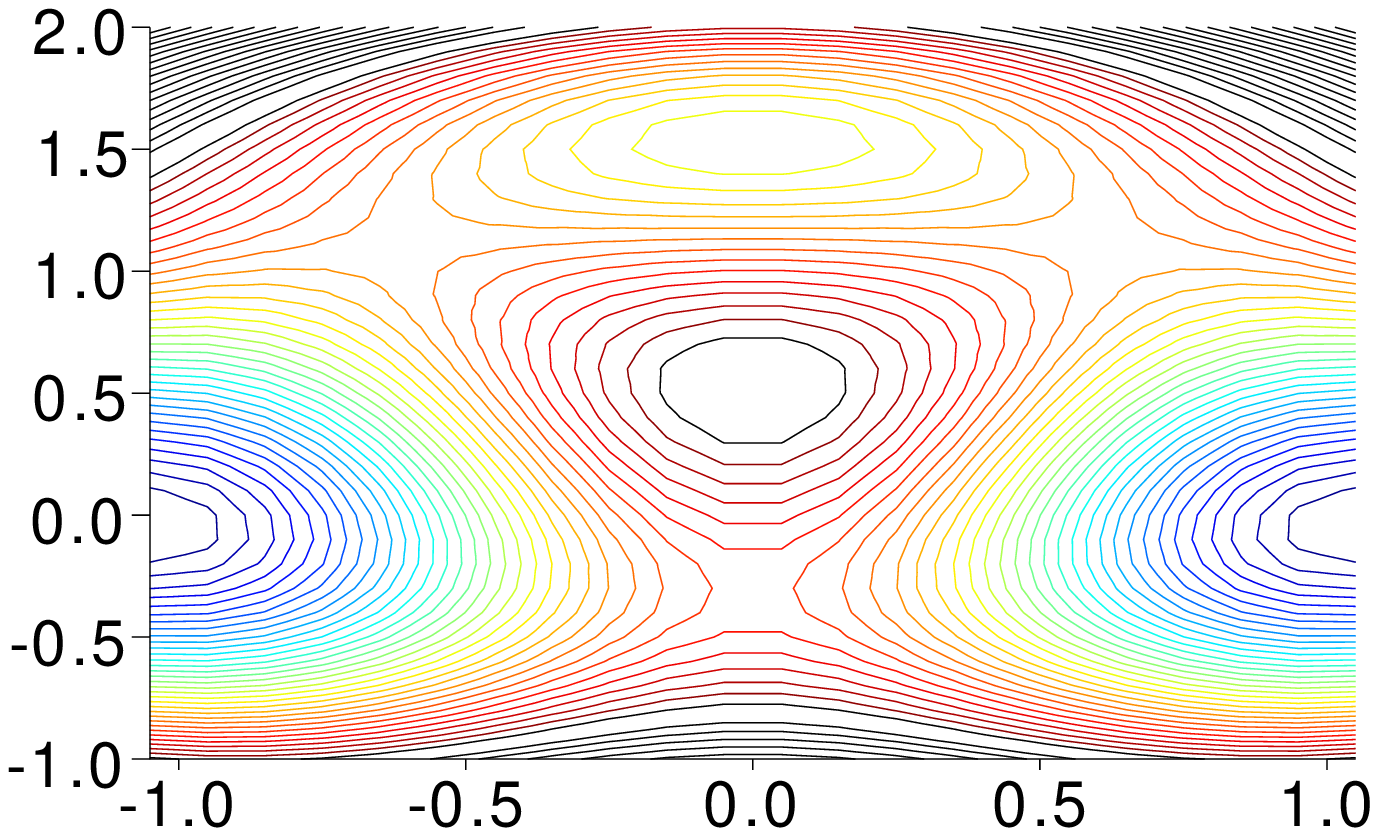}
\includegraphics[width=7.4cm]{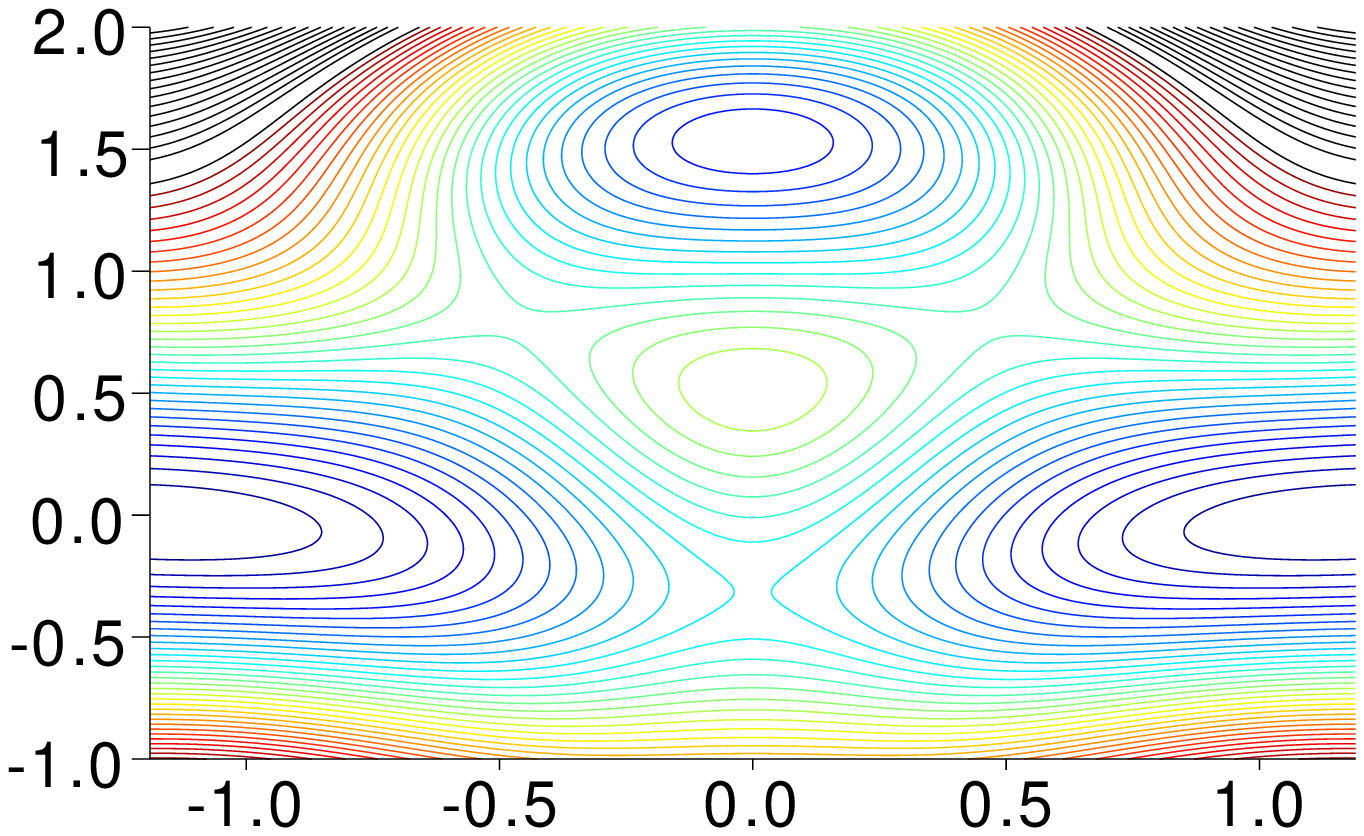}
\caption{\label{fig:contour}
  Left: Level sets of the potential~$U$ defined
  in~\eqref{eq:pot_U}. The minima are located at the positions $x_\pm
  = (\pm 1,0)$, and there are three saddle-points, at the positions $x^{{\rm sd},1}_\pm \simeq (\pm 0.6,1.15)$ and $x^{{\rm sd},2} \simeq (0,-0.3)$. The energy differences of these saddle points with respect to the minimal potential energy are respectively~$\Delta U^1 = 2.2$ and~$\Delta U^2 = 2.7$. 
  Right: Level sets of the biased potential~$U + \beta^{-1} \log \t_\star \propto -\beta^{-1} \log \pi_{\t_\star}$ for $\beta = 20$, $R=1.1$ and $d=55$ ($\t_\star$ being considered as a function with constant values on the strata~$\Xset_\ell$). The position of the saddle point $x^{{\rm sd},2}$ is unaffected, while the saddle points $x^{{\rm sd},1}_\pm$ are shifted to~$(\pm 0.35,0.7)$. The energy differences of the saddle points with respect to the minimal energy are now respectively $\Delta U^{1,{\rm biased}} \simeq 1.65$ and $\Delta U^{2,{\rm biased}} \simeq 1.25$.
}
\end{center}
\end{figure}

A plot of the level sets of the potential~$U$ is presented in
Figure~\ref{fig:contour} (Left). The global minima of the
potential~$U$ are located at the points $x_- = (-1,0)$ and $x_+ =
(1,0)$.  We also provide  a plot of the biased potential associated
with $\pi_{\t_\star}$ (for $\beta = 20$, $R=1.1$ and $d=22$ strata) in Figure~\ref{fig:contour} (Right).

From Laplace's method, the ratio between
the weight of the stratum in the transition region around $x_1 = 0$
and the strata located near the global minima of the potential $U$
(\textit{i.e.} around $x_\pm$) scales like
$\overline{C} \exp(-\beta \mu_0)$ for some positive values
$\overline{C}$ and~$\mu_0$, in the limit $\beta
\to \infty$. In view of~\eqref{eq:theta_star}, we thus expect that the
equivalent of the parameter $\varepsilon$ of
Section~\ref{sec:efficiency} in terms of $\beta$ should be
\begin{equation}
  \label{eq:eps_beta}
  \varepsilon(\beta) = \overline{C} \, \exp(-\beta \mu_0) \eqsp.
\end{equation}
The aim of this section is to check numerically that, assuming this
relation between $\beta$ and $\varepsilon$, the scaling behaviors we obtained in the previous
section on exit times for the very simple toy model with three states
are indeed also observed for a Markovian dynamics with local moves on
the two dimensional potential $U$. Let us now make precise the
dynamics we consider.

The reference (non adaptive) Markov chain $\overline{X}_n$ is obtained by a Metropolis
algorithm, using an isotropic Gaussian proposal with
variance-covariance matrix $\upsilon^2 \, \mathrm{Id}$ where
$\mathrm{Id}$ is the $2 \times 2$ identity matrix. This dynamics is
metastable: for local moves ($\upsilon$ of the order of a fraction of
$\|x_+-x_-\|$, in the following we choose $\upsilon$ in $\{0.025,0.05,0.1,0.2\}$), it takes a lot of  time to go from the left to the
right, or from the right to the left (notice that the potential is
symmetric with respect to the $y$-axis). More precisely, there are two main metastable states: one
located around $x_- = (-1,0)$, and another one around $x_+ = (1,0)$.
These two states are separated by a region of low probability. The
metastability of the dynamics increases with $\beta$ ({\it i.e.} as the temperature
decreases). The larger $\beta$ is, the larger is the ratio between the
weight under $\pi$ of the strata located near the main metastable states and the
weight under $\pi$ of the transition region around $x_1 = 0$, and the more
difficult it is to leave the left metastable state to enter the one on
the right (and conversely).  We compare the reference (non adaptive) Markov chain  to the
associated Wang-Landau dynamics $X_n$. In particular, the proposal kernel used in the
Metropolis algorithm is the same for the Wang-Landau dynamics and for the
reference dynamics. 
As in the previous section, the nonlinear update~\eqref{eq:NL_update}
is used. The step-size sequence is
chosen as in~(\ref{eq:gamma}). The initial weight vector $\t_0$ is
$(1/d,\dots,1/d)$. Notice that the reference dynamics corresponds to
the case when $\gamma_\star = 0$ (no adaption).

\subsection{Expected scalings in the small temperature regime}
\label{sec:expected_scalings}

Average exit times are obtained by performing independent
realizations of the following procedure: initialize the system in the
state $x_-=(-1,0)$, and run the dynamics until the first time index
$\mathcal{N}$ such that $X_{{\mathcal{N}},1} > 1$ (\textit{i.e.} the first component of $X_\mathcal{N}$ is larger than~1). This average exit time is denoted $t_\beta$ for the Wang-Landau dynamics, and $\overline{t_\beta}$ for the reference dynamics.

Before giving the numerical results, let us state the expected scaling behaviors for $\overline{t_\beta}$ and $t_\beta$ in
the limit $\beta \to \infty$, in view of Proposition~\ref{lem:non-adapt},
Proposition~\ref{lem:analytical_adaptive} and~\eqref{eq:eps_beta}. First,
the scaling~\eqref{eq:espTbar13} implies that for the reference dynamics,  under the
relation~\eqref{eq:eps_beta} (in the limit $\beta \to \infty$),
\begin{equation}
  \label{eq:predicted_law_ref}
 \overline{t_\beta} \sim \frac{6}{\overline{C}} \exp(\beta \mu_0).
\end{equation}
Second, for the Wang-Landau dynamics, \eqref{equivt} implies that, under the
relation~\eqref{eq:eps_beta} (in the limit $\beta \to \infty$): for $\alpha \in [1/2,1)$ (and we
will even consider $\alpha \in (0,1)$ below), 
\begin{equation}
  \label{eq:predicted_power_law}
  t_\beta \sim \left(\frac{(1-\alpha)\mu_0}{\gamma_\star} \beta \right)^{1/(1-\alpha)},
\end{equation}
while, for $\alpha = 1$,
\begin{equation}
  \label{eq:predicted_exp_law}
  t_\beta \sim C_{\gamma_\star} \exp\left(\beta \frac{\mu_0}{1+\gamma_\star}\right).
\end{equation}
In practice, the range of values of $\beta$ required to observe the asymptotic
regime $\beta \to \infty$ depends on the values of $\alpha$ and
$\gamma_\star$ (see Figure~\ref{fig:scalings}). 

\subsection{Choice of the numerical parameters} 
\label{sec:choice_param}

For a given value of the inverse temperature~$\beta$, the computed average exit times~$t_\beta$ and $\overline{t_\beta}$ are obtained by averaging over $M$ independent realizations of the process started at $x_-$. We use the Mersenne-Twister random number generator as
implemented in the GSL library. We choose $M$ such that the
relative error on $t_\beta$ or $\overline{t_\beta}$ is less  than a few percents in the worst cases. For computational reasons, $M$ is of the order of a few hundreds for the largest exit times, while $M=10^5$ in the easiest cases. 

The choice of the number of bins is a more delicate matter. We
consider in the sequel $R=1.1$ since we want to observe transitions
from $x_-$ to $x_+$, and decompose the interval $[-R,R]$ into $d$
strata of width $2R/d = \Delta x$. In order to sufficiently refine the
variations of the potential and to produce a not too coarse
free-energy profile, we consider bin widths $\Delta x$ smaller than $0.2$. In order to preserve the locality of the moves, the magnitude
of the random displacements (which are of order~$\upsilon$) is chosen
in order to be comparable to the width $\Delta x$ of one stratum. Therefore, 
from one stratum, the neighboring ones are the most
likely to be visited. This is reminiscent of the dynamics used on the
toy model in the previous section.

\begin{figure}\begin{center}
\includegraphics[width=7.4cm]{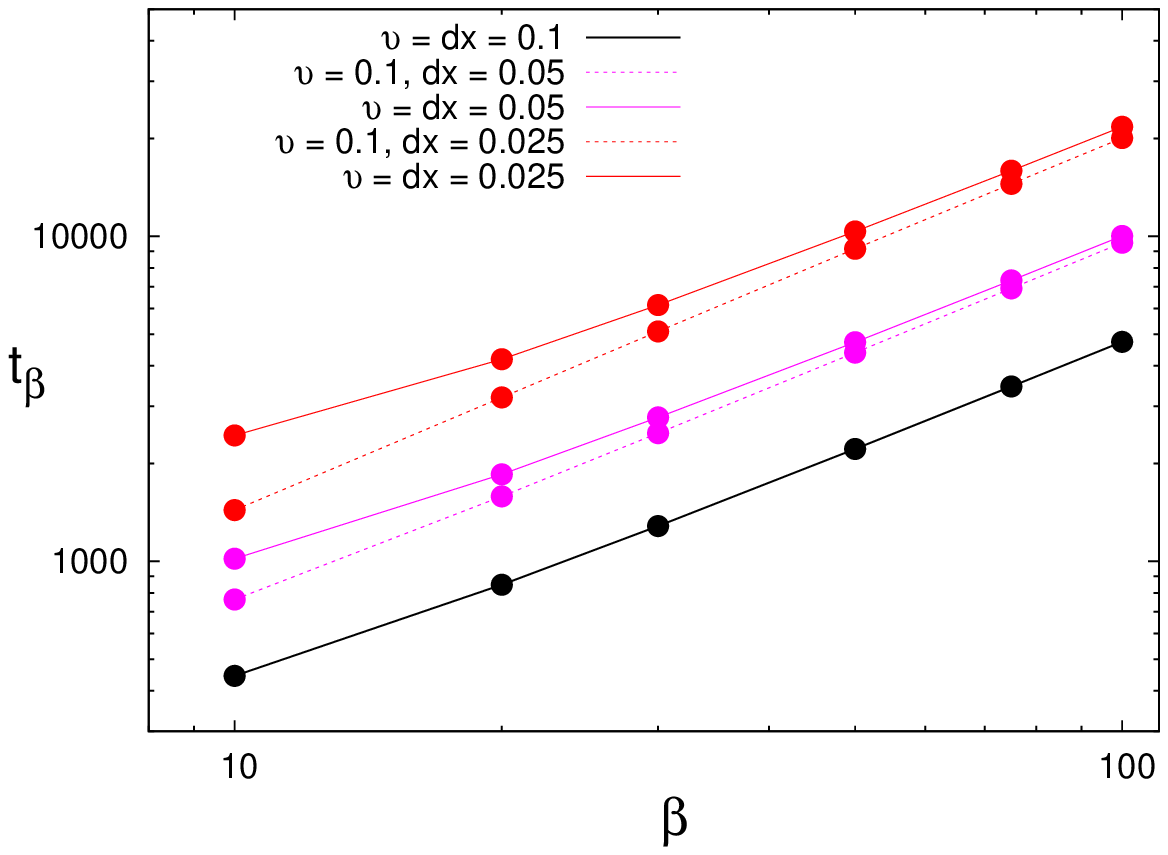}
\includegraphics[width=7.4cm]{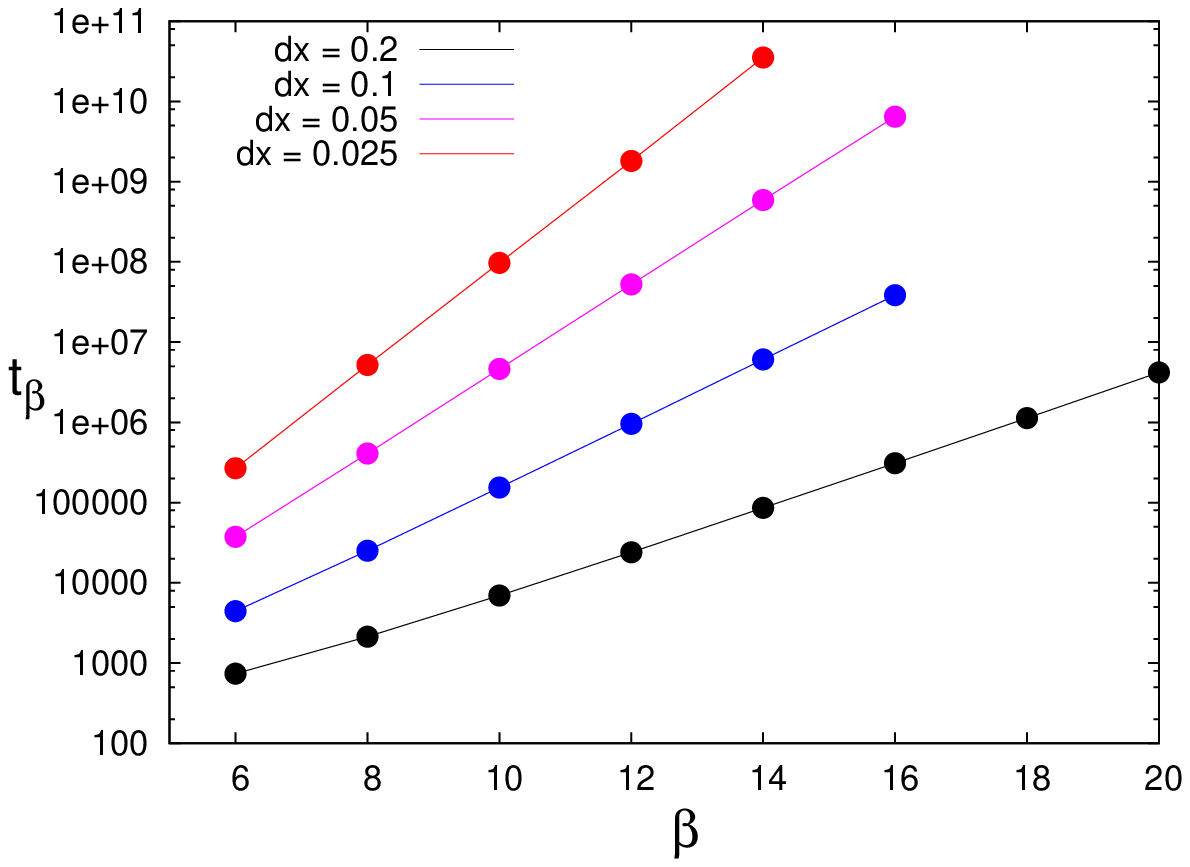}
\caption{\label{fig:dx}
  Left: In the case $\alpha = 0.125$ and $\gamma_\star = 1$  in the stepsize sequence~\eqref{eq:gamma}, the scaling of the average exit times is independent of the number of bins, even if $\upsilon$ is of the order of several $\Delta x$.
  Right: In the case $\alpha = 1$ and $\gamma_\star = 8$, the exponential rate $a$ for the scaling $t_\beta \sim C \exp(a \beta)$ depends on~$\Delta x$ (see Table~\ref{tab:scaling_bin}).
}
\end{center}
\end{figure}

Results on the dependence of the average first exit times $t_\beta$ as a function of $\Delta x$ are presented in Figure~\ref{fig:dx}. The conclusions which can be drawn from these results are the following:
\begin{enumerate}[(i)]
\item when $\alpha = 0.125$ and $\gamma_\star = 1$ (as already hinted at in the beginning of Section~\ref{sec:exit_numerical}, the interest of this case which is not covered by the theoretical analysis of Section~\ref{sec:efficiency} is that the Wang Landau algorithm quickly escapes from metastable states and it is easier to investigate numerically very large values of $\beta$), the average exit time scales in all cases as $t_\beta \sim C \beta^{1/(1-\alpha)}$, as predicted by~\eqref{eq:predicted_power_law}, and only the prefactor depends on the number of bins~$d$. A more precise look at the results shows that the prefactor~$C$ is proportional to~$d$. Note also that the average exit time increases when $\upsilon$ decreases, although this increase is moderate;
\item when $\alpha = 1$ and $\gamma_\star = 8$, the asymptotic behavior depends more dramatically on the number of bins. For all our choices of $\Delta x$, the average exit time scales as $t_\beta \sim C \exp(a \beta)$, as suggested by~\eqref{eq:predicted_exp_law}, but the value $a$ depends on~$\Delta x$. More precisely, the rate $a$ decreases as $\Delta x$ is increased (see the precise results in Table~\ref{tab:scaling_bin}). 
\end{enumerate}
We expect the same conclusions to hold for other values of $\alpha$ and $\gamma_*$, the important distinction being whether $\alpha < 1$ or $\alpha = 1$.

\begin{table}
\begin{center}
\begin{tabular}{|c@{\quad}c|}
\hline
$\Delta x$ & $a$ \\
\hline
0.025 & 1.47 \\  
0.05 & 1.21 \\
0.1 & 0.92 \\
0.2 & 0.63 \\
\hline
\end{tabular}
\end{center}
\caption{\label{tab:scaling_bin}
Fitted value of~$a$ as a function of the bin width $\Delta x = 2R/d$ for the expected scaling relation $t_\beta \sim C
\exp(a \beta)$ corresponding to the data presented in
Figure~\ref{fig:dx} (Right), when $\alpha=1$ and $\gamma_\star=8$.
}
\end{table}

In the sequel (except in Section~\ref{sec:2}), we choose $R = 1.1$ and $d = 22$ in order to have a sufficiently refined free energy profile. Consistently with the above discussion, we set $\upsilon = 0.1$. 

\subsection{Numerical results}
\label{sec:numres}

\begin{figure}[h]
\subfigure[\label{fig:scalings:ref}Reference dynamics (logarithmic
scale on the $y$-axis)]{
\psfrag{beta}{$\beta$}
\psfrag{t}{$\overline{t_\beta}$}
\includegraphics[width=7cm]{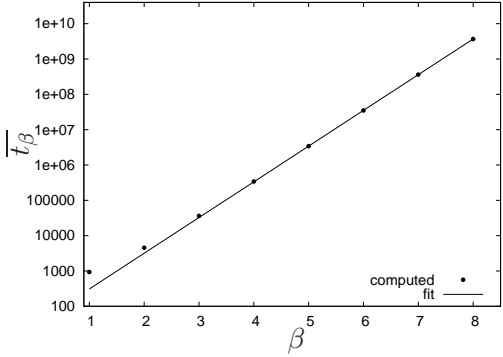}
}
\subfigure[\label{fig:scalings:1_C_2}$\alpha = 1$ and $\gamma_\star = 2$ (logarithmic scale on the $y$-axis)]{
\psfrag{beta}{$\beta$}
\psfrag{t}{${t_\beta}$}
\includegraphics[width=7cm]{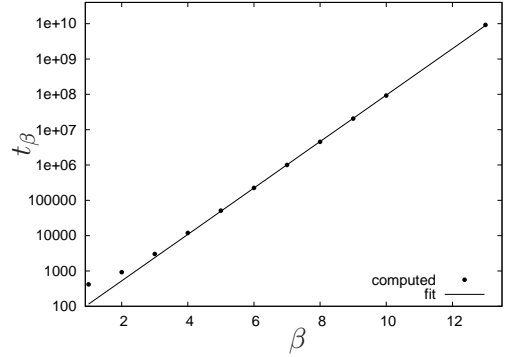}
}
\subfigure[\label{fig:scalings:0.75}$\alpha = 0.75$ (logarithmic scale
on the $x$ and $y$-axis)]{
\psfrag{beta}{$\beta$}
\psfrag{t}{${t_\beta}$}
\includegraphics[width=7cm]{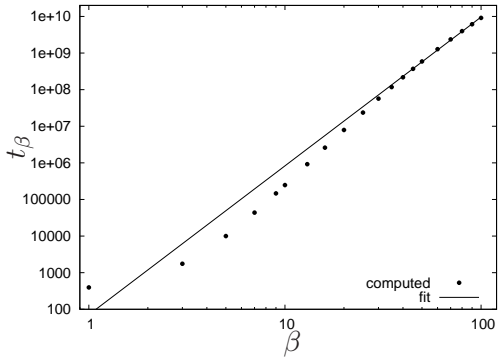}
}
\subfigure[\label{fig:scalings:0.125}$\alpha = 0.125$ (logarithmic
scale on the $x$ and $y$-axis)]{
  \psfrag{beta}{$\beta$}
\psfrag{t}{${t_\beta}$}
\includegraphics[width=7cm]{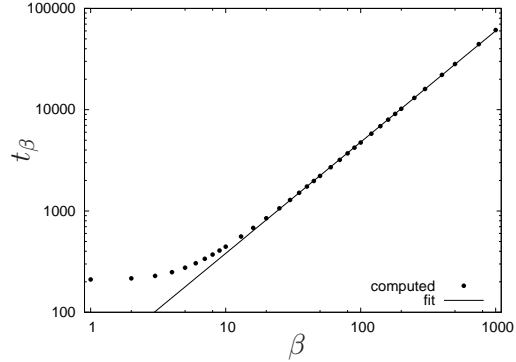}
}
\caption{\label{fig:scalings}
  Average exit time as a function of $\beta$ for various step-size sequences~\eqref{eq:gamma}.
}
\end{figure}

Let us first check that we indeed recover the correct scaling
behavior~\eqref{eq:predicted_law_ref} on the average
exit times for the reference (non adaptive) dynamics. In  Figure~\ref{fig:scalings:ref}, we plot, as a function
of $\beta$, the average exit time $\overline{t_\beta}$ for the non-adaptive
dynamics, using a logarithmic scale on the $y$-axis. The affine fit
is very good, and yields an approximate value for the slope: $\mu_0 \simeq 2.32$.
This value is of the order of the saddle point energy difference $\Delta U^1$ (see the caption 
of Figure~\ref{fig:contour}). 

We then plot $t_\beta$ as a function of $\beta$ in the case $\alpha =
1$ and $\gamma_\star = 2$ in Figure~\ref{fig:scalings:1_C_2}, still
using a logarithmic scale on the $y$-axis. As expected
from~\eqref{eq:predicted_exp_law}, we indeed observe some exponential 
asymptotic behavior $t_\beta \sim C_{\gamma_\star} \exp(\beta \mu_{\gamma_\star})$. 
This is true for other values of $\gamma_\star$. We report the corresponding slopes
$\mu_{\gamma_\star}$ for various values of $\gamma_\star$ in Table~\ref{tab:alpha_one}.
\begin{table}
\begin{center}
\begin{tabular}{|c@{\quad}|c@{\quad}c|}
\hline
$\gamma_\star$ & $\mu_{\gamma_\star}$ & $\mu_0/(1+\gamma_\star)$ \\
\hline
0 & 2.32 & 2.32 \\
1 & 1.74 & 1.16 \\
2 & 1.51 & 0.77 \\
4 & 1.25 & 0.46 \\
8 & 0.92 & 0.26 \\
\hline
\end{tabular}
\end{center}
\caption{\label{tab:alpha_one}
Update with step-sizes $\gamma_n = \gamma_\star/n$ ($\alpha = 1$,
$d=22$ or equivalently, $\Delta x=0.1$). Exponents of the
law $t_\beta \sim C_{\gamma_\star} \exp(\mu_{\gamma_\star}\beta)$ for various values of $\gamma_\star$.}
\end{table}
Although the exponential dependence of $t_\beta$ on~$\beta$ consistent
with~\eqref{eq:predicted_exp_law} is reproduced, the exact dependence
on $\gamma_\star$
of the constant in the exponential predicted by the analytical example
is not exactly observed here since $\mu_{\gamma_\star} \neq \mu_0/(1+\gamma_\star)$. In fact, $\mu_{\gamma_\star}$ is systematically larger than $\mu_0/(1+\gamma_\star)$. This was expected in view of the results presented in Section~\ref{sec:choice_param} (since the exponential rate increases as $\Delta x$ decreases).

We now turn to the case $\alpha \in (0,1)$ where we expect $t_\beta
\sim C_\alpha \beta^{1/(1-\alpha)}$, see~\eqref{eq:predicted_power_law}. 
Note that we also consider the case $\alpha \in (0,1/2)$ which was
not covered by the theoretical analysis of Section~\ref{sec:efficiency}.
To confirm the expected behavior, we plot $t_\beta$ as a function of
$\beta$ in a log-log scale, see
Figure~\ref{fig:scalings:0.75}-\ref{fig:scalings:0.125} for the cases
$\alpha = 0.75$ and $\alpha = 0.125$ respectively. We observe in all
cases a dependence $t_\beta \sim C_\alpha \beta^{\mu_\alpha}$, the
value of the exponent $\mu_\alpha$ being the slope of the affine fit
in the log-log diagram. The estimated exponents are gathered in
Table~\ref{tab:alpha_not_one} for various values of $\alpha$ when
$\gamma_\star = 1$. They compare very well with the value
$1/(1-\alpha)$ predicted from~\eqref{eq:predicted_power_law}. On the
other hand, we were not able to obtain a meaningful dependence of the
prefactor $C_\alpha$ on the parameter~$\alpha$. This is related to the dependence of the prefactor
on the number of bins (see Section~\ref{sec:choice_param}).

In conclusion, these numerical experiments are in very good agreement
with our theoretical findings of Section~\ref{sec:efficiency}.

\begin{table}
\begin{center}
\begin{tabular}{|c@{\quad}|c@{\quad}c|}
\hline
$\alpha$ & $\mu_\alpha$ & $1/(1-\alpha)$\\
\hline
0.125 & 1.11 & 1.14 \\
0.25 & 1.30 & 1.33 \\
0.375 & 1.55 & 1.60 \\
0.5 & 2.02 & 2.00 \\
0.625 & 2.72 & 2.67 \\
0.75 & 4.06 & 4.00 \\
\hline
\end{tabular}
\end{center}
\caption{\label{tab:alpha_not_one}
Update with step-sizes $\gamma_n = n^{-\alpha}$.
Exponents of the scaling law
$t_\beta \sim C_\alpha \beta^{\mu_\alpha}$ for $\alpha \in (0,1)$.
}
\end{table}

\section{Proof of the results presented in Section~\ref{sec:efficiency}}
\label{sec:proof:cvg}

In the following, we denote by $\lfloor x \rfloor$ the integer part of $x \in \mathbb{R}$, namely
the integer such that $\lfloor x \rfloor \leq x < \lfloor x \rfloor +1$. We will
also use the notation $\lceil x \rceil$ for the integer such that $ \lceil x
\rceil - 1 < x \leq \lceil x \rceil$. For $i \neq j \in \{1,2,3\}$, the
time to go from~$i$ to~$j$ for the non-adaptive dynamics is denoted
\begin{equation}\label{eq:Tbar}
\overline{T}_{i \to j}=\min \Big\{n \, : \, \overline{X}_n=j \text{ starting from } \overline{X}_0=i \Big\}.
\end{equation}
A similar definition holds for the time ${T}_{i \to j}$ to go from~$i$ to~$j$ for the Wang-Landau dynamics.

\subsection{Proof of Proposition~\ref{lem:non-adapt}}
\label{sec:proof_non-adapt}

Using the Markov property and decomposing a trajectory from state 1 to state~3 as successive attempts from 1 to 2 back to 1, and eventually a successful transition from 1 to 2 up to 3, it is easy to check that:
\begin{equation}\label{eq:Tbar13}
\overline{T}_{1 \to 3} = \sum_{n=1}^{N} \left( \overline{T}_{1 \to 2}^n + \overline{T}_{2 \to \{1,3\}}^n \right),
\end{equation}
where
\[
N \sim \G\left(\frac12\right), 
\qquad 
\overline{T}_{1 \to 2}^n \sim \G\left(\frac{\varepsilon}{3}\right), 
\qquad 
\overline{T}_{2 \to \{ 1, 3 \}}^n \sim \G\left(\frac23\right),
\]
are independent geometric random variables. The random variable $N$
is the number of jumps from $1$ to $2$ before $3$ is
eventually visited. The random
variables $\overline{T}_{1 \to 2}^n$ (respectively $\overline{T}_{2
  \to \{1,3\}}^n$) are the $n$-th sojourn time in state $1$
(respectively state $2$). Notice that we have used here the fact that
starting from state $2$, the probability to go to state $1$ is equal
to the probability to go to state $3$, which implies that the
parameter of the geometric random variable $N$ is $1/2$.

Let us show that~\eqref{eq:espTbar13} and~\eqref{eq:probTbar13} are easily obtained from~\eqref{eq:Tbar13}.
Indeed, using the fact that for independent geometric random variables $A \sim \G(a)$ and $B_k \sim \G(b)$ (the random variables $B_k$ being i.i.d.),
\[
\sum_{k=1}^{A} B_k \sim \G(ab),
\]
it is easily seen that
$\overline{T}_{1 \to 3} \stackrel{\rm (d)}{=}N_1 + N_2$,
where $N_1$ and $N_2$ are (non-independent) geometric random variables:
\[
N_1 \sim \G\left(\frac{\eps}{6}\right), 
\qquad 
N_2 \sim \G \left(\frac{1}{3}\right).
\]
Therefore $\EE\left(\overline{T}_{1 \to 3}\right)=\frac{6}{\varepsilon}+3$ so that \eqref{eq:espTbar13} holds.
Notice that,  in the limit $\eps \to 0$, we have the following convergences in law:
\[
\eps N_1 \rightarrow \mathscr{E}\left(\frac16\right),
\quad
\eps N_2 \rightarrow 0,
\] 
where $\mathscr{E}(1/6)$ denotes an exponential random variable with
parameter $1/6$. The result~\eqref{eq:probTbar13} is then easily
obtained by the Slutsky theorem.

\subsection{Proof of Proposition~\ref{lem:analytical_adaptive}}
\label{sec:proof_analytical_adaptive}

The heuristic of the proof is the following. In the limit of
small~$\eps$, to go from $1$ to $3$, a typical path first needs to
stay sufficiently long in $1$, in order for a transition to $2$ to be
more likely (when $\theta_n(1)$ becomes sufficiently large). Then, from~$2$, 
the time it takes to go to $3$ is small compared to the time spent to
leave $1$ for the first time. The aim of this proof is to quantify
that by: (i) showing that a transition from $1$ to $2$ in a
well-chosen time is very likely and then (ii) showing that once $2$ is reached, 
the time it remains to go to
$3$ is small compared to the first transition time from~$1$ to~$2$.
The precise result is the following.

\begin{prop}\label{lem:T13}
Consider the Wang-Landau dynamics defined in Section~\ref{sec:def_dyn}. Let us
assume that $\alpha \in [1/2,1]$ and that, if $\alpha=1/2$,
$\gamma_\star < 1$. Then, 
\begin{equation}\label{eq:probT}
\lim_{\eps \to 0} \PP\Big(T_{1 \to 3} \in \big(a(\eps),b(\eps)\big)\Big) =1,
\end{equation}
with
\begin{itemize}
\item for $\alpha \in [1/2,1)$,  
  \[
  a(\eps) = \left(\frac{1-\alpha}{\gamma_\star} \big[
  |\ln \eps| - \beta(\eps) \big]\right) ^{1/(1-\alpha)},
 \qquad
  b(\eps)=C_b |\ln \eps|^{1/(1-\alpha)},
  \]
  where $C_b$ is any constant such that 
  \[
  C_b>\left(\frac{1-\alpha}{\gamma_\star}\right)^{1/(1-\alpha)}
  \]
  and $\beta(\eps)$ is any nonnegative function 
  smaller than $|\ln(\eps)|$ and such that 
  \begin{equation}\label{eq:beta_1}
    \lim_{\eps \to 0} 
    |\ln \eps|^{\alpha/(1-\alpha)} \, \mathrm{e}^{-\beta(\eps)} = 0;
  \end{equation}
\item for $\alpha=1$, 
  \[
  a(\eps)=\eps^{-1/(1+\gamma_\star)} f(\eps), 
  \qquad 
  b(\eps)=\eps^{-1/(1+\gamma_\star)} g(\eps),
  \]
  for any positive functions $f$ and $g$ such that 
  \[
  \lim_{\eps \to 0} f(\eps)=0, 
  \qquad 
  \lim_{\eps \to 0} g(\eps)=\infty.
  \]
\end{itemize}
\end{prop}
In the case $\alpha \in [1/2,1)$, an example of a simple admissible
lower bound is $a(\eps)=C_a |\ln \eps|^{1/(1-\alpha)}$ where $C_a$ is any constant
such that $C_a<\left( \frac{1-\alpha}{\gamma_\star}
\right)^{1/(1-\alpha)}$. In
this case, one should consider $\beta(\eps)=\left(1-\frac{\gamma_\star C_a^{1-\alpha}}{1-\alpha}\right)
|\ln \eps|$ which indeed satisfies~\eqref{eq:beta_1}. Therefore Proposition \ref{lem:T13} implies Proposition~\ref{lem:analytical_adaptive}.

\medskip

Before proving the proposition, let us first introduce some
notation. A convenient rewriting of the Wang-Landau dynamics is:  for all $n \geq 0$, given
$(X_n,\tilde \theta_n)$,
\begin{equation}\label{eq:WL}
\left\{
\begin{aligned}
X_{n+1}& \text{ is sampled according to the kernel } P_{\t_n}(X_n,\cdot) \\
\tilde\theta_{n+1}(i)&=\tilde\theta_n(i) ( 1+ \gamma_{n+1} \un_{X_{n+1}=i})
\end{aligned}
\right.
\end{equation}
where $\tilde\theta_0=(1,1,1)$, $P_\t$ is defined by~\eqref{eq:Ptheta} and the normalized weights $\theta_n$
associated with the unnormalized weights $\tilde{\theta}_n$
 are
\[
\theta_n(i)=\left(\sum_{j=1}^3 \tilde\theta_n(j)\right)^{-1} \tilde\theta_n(i). 
\]
The updating rule in~\eqref{eq:WL} is exactly the standard update~\eqref{eq:NL_update}.

A crucial role will be played
by the time the dynamics needs to first reach $2$:
\begin{equation}\label{eqdef:T120}
{T}^0_{1 \to 2}=\min\Big\{n \, : \  {X}_n=2 \text{ starting from } {X}_0=1 \Big\}.
\end{equation}
The probability to go from state~$1$ to state~$2$ in exactly~$n$ moves is
\begin{equation}\label{eq:probT12}
\PP(T^0_{1\to 2} = n) = p_{11}^0 \dots p_{11}^{n-2} p_{12}^{n-1},
\end{equation}
with
\[
p_{11}^m = 1-\frac{1}{3}
\left(\eps \Xi_m \wedge 1 \right), 
\qquad
p_{12}^m = 1-p_{11}^m = \frac{1}{3}
\left(\eps \Xi_m \wedge 1 \right),
\]
where
\begin{equation}\label{eq:Xi}
\Xi_m=\prod_{k=1}^{m} (1+\gamma_k)
\end{equation}
with the convention $\Xi_0=1$.
The first $n-1$ factors in~\eqref{eq:probT12} correspond to staying in state~1 (with the appropriate update
of the weights), and the last one corresponds to the transition from state~1 to state~2.
An important inequality, which will be used below, is
$p_{12}^m \leq p_{12}^n$ (and thus $p_{11}^m \geq p_{11}^n$) for $m \leq n$: when the system is stuck in state 1, as time goes, the probability to go to state~$2$ increases. 

\begin{figure}[htbp]
\begin{center}
\psfrag{a}{$T^0_{1\to 2}$}
\psfrag{b}{$\tau^0_{1\to 2}$}
\psfrag{c}{$\tau^1_{2\to 1}$}
\psfrag{d}{$T^1_{1\to 2}$}
\psfrag{e}{$\tau^1_{1\to2}$}
\psfrag{h}{$\tau^{N_{2 \to 1}}_{2\to 1}$}
\psfrag{g}{$T^{N_{2 \to 1}}_{1\to 2}$}
\psfrag{f}{$\tau^{N_{2\to 1}}_{1\to2}$}
\psfrag{i}{$T_{1\to 3}$}
\includegraphics[width=\textwidth]{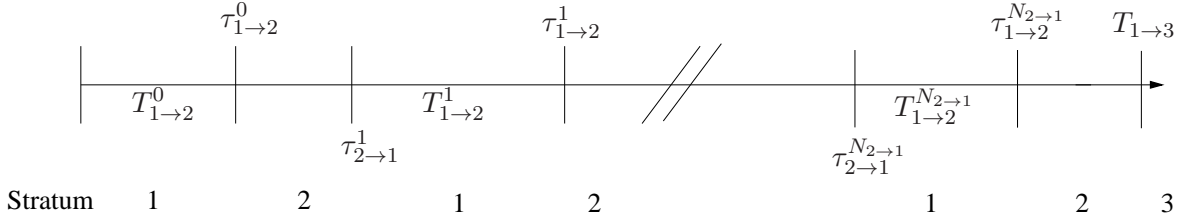}
\caption{\label{fig:exit_time}
  Schematic representation of the successive passage times and exit times out of~1.
}
\end{center}
\end{figure}

Estimates on the exit time $T_{1\to 3}$ are based on the following equality:
\begin{equation}\label{eq:T13equality}
T_{1\to 3}=T^0_{1\to2}+\sum_{i=1}^{N_{2\to1}}T^i_{1\to2}+N_2
\end{equation}
where $N_2=\sum_{n=0}^{T_{1 \to 3}} 1_{\{X_n=2\}} $ is the time the chain spends in $2$ before going to
$3$, $N_{2 \to 1}$ is the number of jumps from $2$ back to $1$ before
going to $3$ and $T^i_{1\to2}$ is the time it takes to leave $1$ at
the $i$-th return to the state $1$ from $2$. Notice that
\[
N_{2 \to 1} \leq N_2.
\]
To make these quantities more precise, let us introduce the successive passage times: 
for $i\geq 1$,
\begin{equation}\label{eq:taui21}
\tau^i_{2\to 1}=\inf \left\{n > \tau^{i-1}_{1 \to 2}, \, X_n=1\right\},
\end{equation}
with, by convention $\tau^{0}_{1 \to 2}=T^0_{1 \to 2}$ and,
\[
\tau^i_{1\to 2}=\inf \left\{n > \tau^{i}_{2 \to 1}, \, X_n=2\right\}.
\]
Note that $T^i_{1 \to 2}=\tau^{i}_{1 \to 2}-\tau^i_{2 \to 1}$. We
refer to Figure~\ref{fig:exit_time} for a schematic representation of
all these times.

Let us first state a simple result concerning $N_2$ which is based on the fact that before visiting the state $3$,
$\tilde \theta_n(3)=1$ remains unchanged while $\tilde \theta_n(2) \geq 1$. This means that for $n \leq T_{1 \to 3}$, 
\[
P_{\theta_n}(2,3)=\frac{1}{3} \left(\frac{\tilde
\theta_n(2)}{\eps\tilde
\theta_n(3)}  \wedge 1 \right)=\frac{1}{3} \left(\frac{\tilde
\theta_n(2)}{\eps}  \wedge 1 \right)=\frac{1}{3},
\] where we have used the inequality $\eps < 1$. At each time the system is in state $2$, it stays in state $2$ or goes to state $1$ at the next time with probability $2/3$. This gives the intuition of the following result, the formal proof of which is postponed to Section \ref{sec:tecres}.  
\begin{lemma}\label{lem:N2}
The random variable $N_2$ is geometric with parameter $1/3$: for all $n \geq 0$,
$$\PP(N_2 \geq n) = \left(\frac23\right)^n.$$
\end{lemma}

Thus, in~\eqref{eq:T13equality}, the last term plays no role in
the limit $\eps \to 0$. We show below that this is also true for
the second term: the main role is played by $T^0_{1 \to 2}$. This is
why we first need to precisely estimate the time $T^0_{1 \to 2}$.
This can be done for any $\alpha \in (0,1]$ (and not only in
$[1/2,1]$), and without any restriction on $\gamma_\star$.

\begin{lemma}
\label{lem:T120}
Fix $\alpha \in (0,1]$. For $\alpha \in [1/2,1]$, let $a$ be the function defined in
Proposition~\ref{lem:T13}. For $\alpha \in (0,1/2)$, let the function
$a$ be defined in the same way as for $\alpha \in [1/2,1)$. Then, 
\begin{equation}\label{eq:T120}
\lim_{\eps \to 0} \PP\Big(T^0_{1 \to 2} \in \left(a(\eps),\tilde b(\eps)\right) \Big) =1,
\end{equation}
where
\begin{itemize}
\item if $\alpha \in (0,1)$, $\tilde b(\eps)=C_{\tilde b} |\ln
  \eps|^{1/(1-\alpha)}$  where $C_{\tilde b}$ is any constant such that $\dps C_{\tilde b}>\left(\frac{1-\alpha}{\gamma_\star}\right)^{1/(1-\alpha)}$; 
\item if $\alpha=1$, $\tilde b(\eps)=\tilde g(\eps) \eps^{-1/(1+\gamma_\star)}$ for any positive
  function $\tilde{g}$ such that $\dps \lim_{\eps \to 0} \tilde g(\eps)=\infty$.
\end{itemize}
\end{lemma}

 The proof of Lemma~\ref{lem:T120} can be read in 
 Section~\ref{sec:tecres}.

\begin{rem}
To guess the correct scaling for $T^0_{1 \to 2}$, one may consider
the typical time $n(\eps)$ for which $\PP(T^0_{1 \to 2}\leq n(\eps))=1-\prod_{k=0}^{n(\eps)-1} \left(1 - \frac13 \left(\eps
    \Xi_{k} \wedge 1 \right)\right)$ has a positive limit when $\eps$
goes to $0$. Using an expansion when $\eps$ goes to $0$, assuming that $\eps
\Xi_{n(\eps)}$ goes to zero,
we obtain that $n(\eps)$ satisfies
$\sum_{k=0}^{n(\eps)-1} \Xi_{k} \sim \frac{C}{\eps}$ for some constant
$C>0$. A guess for the scaling of the time $T^0_{1 \to 2}$ is thus
$n(\eps) = \arg\min_n \left\{ \sum_{k=0}^{n-1} \Xi_k \geq
  \frac{1}{\eps} \right\}$ (obtained by choosing $C=1$).
Using Lemma~\ref{lem:XiLB} below, this yields various asymptotic behaviors for
$n(\eps)$ depending on the values of $\alpha$ and $\gamma_\star$
in~\eqref{eq:gamma}: 
\begin{itemize}
\item When $\alpha \in (0,1)$, from~\eqref{eq:equivalent_lnXi_alpha<1}, we
   obtain that $n(\eps) \sim \left(\frac{1-\alpha}{\gamma_\star}\right)^{1/(1-\alpha)} |\ln
\eps|^{1/(1-\alpha)}$.
\item When $\alpha=1$, from~\eqref{eq:prod_cas1}, we obtain 
$n(\eps) \sim
\Gamma(2+\gamma_\star)^{1/(1+\gamma_\star)}
\, \eps^{-1/(1+\gamma_\star)}$.
\end{itemize}
This motivates the scaling for $T^0_{1 \to 2}$.
\end{rem}

 We are now in position to prove
 Proposition~\ref{lem:T13}.
\begin{adem}[ of Proposition \ref{lem:T13}]
Let $a$ and $b$ satisfy the assumptions of Proposition \ref{lem:T13}.
Since $T_{1\to 3} \geq T^0_{1\to2}+1$, the lower bound in Proposition~\ref{lem:T13} (i.e. the fact that $\lim_{\eps\to 0}\PP(T_{1 \to 3}\leq a(\eps))=0$) immediately follows from Lemma \ref{lem:T120}. The upper bound
requires some more work. We choose $\tilde b(\eps)$
satisfying the assumptions of Lemma~\ref{lem:T120} with $\left(\frac{1-\alpha}{\gamma_\star}\right)^{1/(1-\alpha)}<C_{\tilde b}<C_b$ if $\alpha\in [1/2,1)$ and $\tilde{g}<g$ if $\alpha=1$. In particular, $\tilde b<b$. Let us also
introduce 
a positive function $\Delta(\eps)$ going to infinity as $\eps\to 0$,
that will be specified later on. Then, using~\eqref{eq:T13equality},
\begin{align}
\PP(T_{1 \to 3} \geq b(\eps)) 
& \leq \PP\Big( T^0_{1 \to 2} \notin \left(a(\eps),\tilde b(\eps)\right)\Big) 
+ \PP\Big(N_2 \geq \Delta(\eps)\Big) \nonumber \\
& \quad + \PP\Big(T^0_{1 \to 2} \in \left(a(\eps),\tilde b(\eps)\right), \, N_2 \leq \Delta(\eps), \, T_{1 \to 3} \geq b(\eps) \Big)  \nonumber \\
& \leq \PP\Big( T^0_{1 \to 2} \notin \left(a(\eps),\tilde b(\eps)\right)\Big) + \PP(N_2 \geq \Delta(\eps)\Big)\label{eq:majT13}
\\
& \quad + \PP\left(T^0_{1 \to 2} \in \left(a(\eps),\tilde b(\eps)\right), \, N_2 \leq \Delta(\eps),
\, \sum_{i=1}^{N_{2\to1}}T^i_{1\to2} \geq b(\eps) -
\Delta(\eps) - \tilde b(\eps) \right). \nonumber
\end{align}
The first term in the right-hand side goes to zero as $\eps \to 0$ by Lemma \ref{lem:T120}. Since $\Delta(\eps)$ tends to
$\infty$ when $\eps$ goes to zero, the second term goes to zero
by Lemma~\ref{lem:N2}. Concerning the third term, the idea is the following: we would like to choose 
$\tilde b$ and $\Delta$ such that, on the event $T^0_{1 \to 2} \in \left(
a(\eps),\tilde b(\eps)\right)$ and  $N_2 \leq \Delta(\eps)$, the times
$T^i_{1\to2}$ can be simply controlled using the fact that the state
$1$ has already been visited for a long time (namely $T^0_{1 \to 2} >
a(\eps)$ and therefore $\tilde\theta_{T^0_{1 \to 2}}(1)$ is large) and the state $2$ is not visited many times (this corresponds to $N_2
\leq \Delta(\eps)$ so that $\tilde\theta_n(2)$ remains small).
This idea will be quantified in Lemma \ref{lem:Delta} in Section \ref{sec:tecres} from which we will deduce :
\begin{lemma}\label{lem:t3}
Assume that $\Delta(\eps)={\rm O}(a^\alpha(\eps))$ as $\eps\to 0$. Then,
there exist constants $C,C',\bar{\varepsilon}>0$ such that for $\varepsilon\in (0,\bar{\varepsilon})$,
\begin{align*}
   & \PP\left(T^0_{1 \to 2} \in \left(a(\eps),\tilde b(\eps)\right), \, N_2 \leq\Delta(\eps),
   \, \sum_{i=1}^{N_{2\to1}}T^i_{1\to2} \geq b(\eps) -
\Delta(\eps) - \tilde b(\eps) \right)\\
 &\qquad \qquad \qquad \qquad\leq \begin{cases}
   \dps \Delta(\eps)\exp\left(-C \frac{1}{\Delta(\eps)} |\ln \eps|^{2- \gamma_\star^2}
  \exp \left( - \beta(\eps)\right)\right)\quad\mathrm{ if } \ \alpha=1/2 \\\Delta(\dps \eps)\exp\left(-C \frac{1}{\Delta(\eps)} |\ln \eps|^{1/(1-\alpha)}
  \exp \left( - \beta(\eps)\right)\right)\quad\mathrm{ if } \ \alpha\in (1/2,1)\\\dps C'\Delta(\eps) \exp \left( - C \frac{ g(\eps) - \tilde g(\eps)}{\Delta(\eps)} 
f(\eps)^{\gamma_\star}\right)\quad\mathrm{ if } \ \alpha=1.
\end{cases}
\end{align*}
\end{lemma}

Let us first conclude in the case $\alpha \in [1/2,1)$.  
We may choose $\Delta$ satisfying $\Delta(\eps)={\rm O}(a^\alpha(\eps))$ and going to infinity as slowly as
needed. Then, for the upper-bound of the third term of the right-hand-side of \eqref{eq:majT13} given by Lemma \ref{lem:t3} to vanish as $\eps\to 0$, it is enough that
\begin{align*}
   & \lim_{\eps \to 0} |\ln \eps|^{1/(1-\alpha)}
  \exp \left( - \beta(\eps)\right)=+\infty\mbox{ if }\alpha \in (1/2,1),\\
\qquad & \lim_{\eps \to 0}|\ln \eps|^{2-\gamma_\star^2}
  \exp \left( - \beta(\eps)\right)=+\infty\mbox{ if }\alpha=1/2.
\end{align*}
For this limit to hold, it is always possible to decrease $\beta$ as long as \eqref{eq:beta_1} holds i.e.  $$\lim_{\eps \to 0} 
    |\ln \eps|^{\alpha/(1-\alpha)} \, \mathrm{e}^{-\beta(\eps)} = 0$$
since the smaller $\beta$, the larger $a$ and the stronger the conclusion \eqref{eq:probT} of Proposition \ref{lem:T13}. 
This is possible without restriction when $\alpha\in (1/2,1)$ and if and only if $2-\gamma_\star^2 > 1$ when $\alpha=1/2$.

Let us now suppose that $\alpha=1$. Up to increasing $f$, which makes
the conclusion of Proposition~\ref{lem:T13} stronger, while preserving
$\lim_{\eps\to 0}f(\eps)=0$, one may assume that $g(\eps) f(\eps)^{\gamma_\star}$ goes to infinity as $\eps$ goes to
zero. In addition, it is always possible to choose $\tilde g \leq g$ such
that $(g(\eps) - \tilde g(\eps))  f(\eps)^{\gamma_\star}$ goes to infinity as $\eps$ goes to
zero. Then one can choose $\Delta$ which grows sufficiently slowly at
infinity so that $\Delta(\eps) \exp \left( - C \frac{ g(\eps) - \tilde g(\eps)}{\Delta(\eps)} 
f(\eps)^{\gamma_\star}\right)$ tends to $0$.
\end{adem}

\subsection{Proofs of the technical Lemmas}
\label{sec:tecres}
\begin{adem}[~of Lemma \ref{lem:N2}]
Let $({\mathcal F}_{n})_{n \geq 0}$ denote the filtration generated by
the Markov chain $\{(X_n,\tilde\theta_n), n \geq 0\}$. Let us also introduce the successive visit times of state $2$. For $i\geq 1$ let $\eta_i=\inf\{n>\eta_{i-1}:X_n=2\}$ with convention $\eta_0=0$.
For $n\in\N^*$, one has $\{N_2\geq n+1\}=\cap_{k=1}^n\{X_{\eta_k+1}\in\{1,2\}\}$. Therefore
\begin{align*}
   \PP(N_{2}\geq n+1)&=\EE\left(\EE\left(\prod_{k=1}^n1_{\{X_{\eta_k+1}\in\{1,2\}\}}|{\mathcal F}_{\eta_n}\right)\right)=\EE\left(\prod_{k=1}^{n-1}1_{\{X_{\eta_k+1}\in\{1,2\}\}}\PP(X_{\eta_n+1}\in\{1,2\}|{\cal F}_{\eta_n})\right)\\&=\EE\left(\prod_{k=1}^{n-1}1_{\{X_{\eta_k+1}\in\{1,2\}\}}(1-P_{\theta_{\eta_n}}(2,3))\right),
\end{align*}
where we used that the event
$\cap_{k=1}^{n-1}\{X_{\eta_k+1}\in\{1,2\}\}$ is ${\cal
  F}_{\eta_n}$-measurable for the second equality and the strong
Markov property for the chain $\{(X_l,\tilde \theta_l),l\geq 0\}$ for the last equality. On $\cap_{k=1}^{n-1}\{X_{\eta_k+1}\in\{1,2\}\}$, the sequence $\{X_l,l\geq 0\}$ has not visited state $3$ before the stopping time $\eta_n$, which implies $P_{\theta_{\eta_n}}(2,3)=\frac{1}{3}$.
Hence $$\PP(N_{2}\geq n+1)=\frac{2}{3}\EE\left(\prod_{k=1}^{n-1}1_{\{X_{\eta_k+1}\in\{1,2\}\}}\right)=\frac{2}{3}\PP(N_{2}\geq n)$$
and one concludes by induction on $n$.
\end{adem}
To prove Lemmas \ref{lem:T120} and \ref{lem:t3}, we need the following estimations on $\Xi_n$.\begin{lemma}\label{lem:XiLB}
For $\alpha=1$,
\begin{equation}
\label{eq:prod_cas1}
\Xi_n\sim \frac{n^{\gamma_\star}}{\Gamma(1+{\gamma_\star})}\mbox{ as }n\to\infty. 
\end{equation}
For $\alpha\in (0,1)$, 
\begin{align}
&\ln (\Xi_n)\sim \frac{\gamma_\star}{1-\alpha}n^{1-\alpha}\mbox{ as
}n\to\infty \label{eq:equivalent_lnXi_alpha<1}\\
&\forall n,\;\Xi_{n}\leq\exp \left( \frac{\gamma_\star}{1-\alpha}
  n^{1-\alpha}\right)\label{eq:XiUB} 
\end{align}
and there exists a constant $C > 0$ independent of $n$ such that
\begin{equation}
 \forall n,\;\Xi_n \geq  \begin{cases}C\exp \left( 2 \gamma_\star  \sqrt{n}- \frac{\gamma_\star^2}{2} \ln n \right)\mbox{ for }\alpha=1/2\\
C\exp \left( \frac{\gamma_\star}{1-\alpha} n^{1-\alpha}\right)\mbox{ for }\alpha \in (1/2,1)\end{cases}.
\label{eq:XiLB} \end{equation}\end{lemma}
\begin{adem}[~of Lemma \ref{lem:XiLB}]
In the case $\alpha=1$, using the Stirling formula, we have
\begin{equation*}
\Xi_n=\prod_{k=1}^n (1+\gamma_k) =\prod_{k=1}^n \left(1+\frac{{\gamma_\star}}{k}\right) 
=\frac{\Gamma (n+1+{\gamma_\star})}{\Gamma(1+{\gamma_\star}) \Gamma(n+1)} 
\sim \frac{n^{\gamma_\star}}{\Gamma(1+{\gamma_\star})}, 
\end{equation*}
which is~\eqref{eq:prod_cas1}. Now, for $\alpha \in (0,1)$, as $n\to\infty$,
\begin{equation*}
\ln (\Xi_n)= \ln \left( \prod_{k=1}^n (1+\gamma_k) \right) \sim \gamma_\star \sum_{k=1}^n k^{-\alpha} \sim \frac{\gamma_\star}{1-\alpha}n^{1-\alpha}.
\end{equation*}
Moreover, 
\begin{align*}
\ln ( \Xi_{n} )
\leq \gamma_\star \sum_{k=1}^{n} k^{-\alpha} 
\leq \gamma_\star \sum_{k=1}^{n} \int_{k-1}^k x^{-\alpha} \, dx 
= \frac{\gamma_\star}{1-\alpha} n^{1-\alpha}. 
\end{align*}
To prove \eqref{eq:XiLB}, we start from the lower bound
\begin{align*}
\ln (\Xi_n) \geq \sum_{k=1}^n \gamma_k - \frac{1}{2} \sum_{k=1}^n \gamma_k^2.
\end{align*}
For $\alpha \in (0,1)$, 
\begin{align*}
\sum_{k=1}^n \gamma_k \geq \gamma_\star \sum_{k=1}^n \int_{k}^{k+1} x^{-\alpha} \, dx 
= \frac{\gamma_\star}{1-\alpha} \left( (n+1)^{1-\alpha} - 1 \right)
\geq \frac{\gamma_\star}{1-\alpha} \left( n^{1-\alpha} - 1 \right),
\end{align*}
so that
\begin{align*}
\ln (\Xi_n) \geq \frac{\gamma_\star}{1-\alpha} \left( n^{1-\alpha} - 1 \right)  - \frac{1}{2} \sum_{k=1}^n \gamma_k^2.
\end{align*}
We now distinguish between two cases.
For $\alpha \in (1/2,1)$, 
\begin{align*}
\sum_{k=1}^n \gamma_k^2
&= \gamma_\star^2 \sum_{k=1}^n k^{-2\alpha} 
\leq \gamma_\star^2 + \gamma_\star^2 \sum_{k=2}^n \int_{k-1}^k x^{-2\alpha} \, dx 
= \gamma_\star^2 + \frac{\gamma_\star^2}{2 \alpha-1}\left( 1 -  n^{1-2\alpha}\right) 
\leq  \frac{2 \gamma_\star^2\alpha}{2 \alpha-1}. 
\end{align*}
Therefore, for $n \geq 1$,
\begin{align*}
\ln (\Xi_n) \geq \frac{\gamma_\star}{1-\alpha} \left( n^{1-\alpha} - 1 \right) - \frac{\gamma_\star^2\alpha }{2 \alpha-1},
\end{align*}
which gives the expected result. For $\alpha=1/2$, 
\begin{align*}
\sum_{k=1}^n \gamma_k^2
\leq \gamma_\star^2 + \gamma_\star^2 \sum_{k=2}^n \int_{k-1}^k x^{-1} \, dx
= \gamma_\star^2 \left( 1 + \ln n\right),
\end{align*}
so that, for $n\geq 1$,
\begin{equation*}
\ln ( \Xi_n ) \geq 2\gamma_\star \left( \sqrt{n} - 1 \right) -  \frac{\gamma_\star^2}{2} \left( 1 + \ln n\right),
\end{equation*}
which also gives the claimed result.
\end{adem}\begin{adem}[~of Lemma \ref{lem:T120}]
Let us first deal with $\alpha \in (0,1)$. 
We start by the lower bound on $T^0_{1\to 2}$. Let $a$ be of
the form $a(\eps) = \left(\frac{1-\alpha}{\gamma_\star} \left(
|\ln \eps| - \beta(\eps) \right)
\right) ^{1/(1-\alpha)}$ for any non-negative function $\beta(\eps)$
smaller than $|\ln(\eps)|$ and satisfying~\eqref{eq:beta_1}. 
By \eqref{eq:probT12}, 
\begin{align*}
\ln \Big( \PP\big\{T^0_{1\to 2} >  a(\eps)\big\} \Big) &= \ln \left(
  \prod_{k=0}^{\lfloor a(\eps)\rfloor} p_{11}^k \right) 
  = \sum_{k=0}^{\lfloor a(\eps)\rfloor} \ln \left(1-\frac{1}{3}
\left(\eps \Xi_k \wedge 1 \right) \right)
\geq - \frac{C_0}{3} \sum_{k=0}^{\lfloor a(\eps)\rfloor}
\left(\eps \Xi_k \wedge 1 \right) \\
& \geq - \frac{C_0 \eps}{3} \sum_{k=0}^{\lfloor a(\eps)\rfloor} \Xi_k,
\end{align*}
where we have used that, by concavity of the function $\ln$, $\ln (1-x) \geq -
C_0 x$ for $x \in (0,1/3)$ with $C_0 = - 3 \ln(2/3) > 0$. 
Now, by \eqref{eq:XiUB},
\begin{align*}
\sum_{k=0}^{n} \Xi_k
&\leq \sum_{k=0}^{n} \exp \left( \frac{\gamma_\star}{1-\alpha}
  k^{1-\alpha} \right)
  \leq \sum_{k=0}^{n} \int_k^{k+1} \exp \left( \frac{\gamma_\star}{1-\alpha}
  x^{1-\alpha} \right)\, dx \\
&  = \int_0^{n+1} \exp \left( \frac{\gamma_\star}{1-\alpha}
  x^{1-\alpha} \right) \, dx 
\leq \frac{(n+1)^\alpha}{\gamma_\star} \int_0^{n+1} \gamma_\star x^{-\alpha} \exp \left( \frac{\gamma_\star}{1-\alpha}  x^{1-\alpha} \right) \, dx \\
& \leq \frac{1}{\gamma_\star}  (n+1)^\alpha \exp \left( \frac{\gamma_\star}{1-\alpha}(n+1)^{1-\alpha} \right).
\end{align*}
Hence, using the inequality $(x+y)^{\delta} \leq x^{\delta} +
y^{\delta}$ for any $(x,y) \in \mathbb{R}_+^2$ and $\delta \in (0,1)$,
\begin{align*}
\ln \Big( \PP\big\{T^0_{1\to 2} > a(\eps)\big\} \Big)
& \geq - \frac{C_0 \eps}{3 \gamma_\star}\,  (a(\eps)+1)^\alpha \exp \left( \frac{\gamma_\star}{1-\alpha}
  (a(\eps)+1)^{1-\alpha} \right)\\
& \geq - C_1 \eps  a(\eps) ^\alpha \exp \left( \frac{\gamma_\star}{1-\alpha}
  a(\eps)^{1-\alpha} \right),
\end{align*}
where $C_1$ is a constant independent of $\eps$. Therefore,
\begin{align}
\ln \Big( \PP\big\{T^0_{1\to 2} > a(\eps)\big\} \Big)
& \geq - C_1 \eps  \left(\frac{1-\alpha}{\gamma_\star} \left(
|\ln \eps| - \beta(\eps) \right)
\right) ^{\alpha/(1-\alpha)} \exp \left(|\ln \eps| - \beta(\eps)
\right) \nonumber \\
& =  - C_2  \left(
|\ln \eps| - \beta(\eps) \right)^{\alpha/(1-\alpha)} \exp \left(- \beta(\eps)
\right), \label{eq:C2}
\end{align}
where $C_2 = C_1 \left(\frac{1-\alpha}{\gamma_\star}\right)^{\alpha/(1-\alpha)}$ is a constant independent of $\eps$.
Thus, under the assumption~\eqref{eq:beta_1}, we indeed obtain that $\lim_{\eps \to 0} \PP( T^0_{1 \to 2} \le
a(\eps)) = 0$.

We now turn to an estimate of an upper bound for $T^0_{1\to2}$. Let us
introduce a function $\tilde b(\eps)=C_{\tilde b} |\ln
  \eps|^{1/(1-\alpha)}$  where $C_{\tilde b}$ is any constant such that  $C_{\tilde b}>\left(\frac{1-\alpha}{\gamma_\star}\right)^{1/(1-\alpha)}$.
We also define an intermediate time $\tilde{n}(\eps) \leq \tilde{b}(\eps)$ 
such that $p_{12}^{\tilde{n}(\eps)}=1/3$, which equivalently writes
\[
\frac13 \left( \eps \Xi_{\tilde{n}(\eps)} \wedge 1
\right)=\frac13\;\;\mbox{ i.e. }\;\;\eps \Xi_{\tilde{n}(\eps)} \geq 1.
\]
We choose
\[
\tilde{n}(\eps)= \left\lceil \tilde{C} |\ln \eps|^{1/(1-\alpha)} \right\rceil, 
\qquad 
\left(\frac{1-\alpha}{\gamma_\star}\right)^{1/(1-\alpha)}<\tilde{C}<C_{\tilde b}. 
\]
In view of~\eqref{eq:equivalent_lnXi_alpha<1}, and since $\tilde{C}^{\alpha - 1} < \frac{\gamma_\star}{1-\alpha}$, it holds $\Xi_n \geq \exp( \tilde{C}^{\alpha - 1} n^{1-\alpha})$ for $n$ large enough. Thus, for $\eps$ small enough, we obtain
\[
\Xi_{\tilde{n}(\eps)} \geq \exp\left( \tilde{C}^{\alpha - 1} \left\lceil \tilde{C} |\ln \eps|^{1/(1-\alpha)}\right\rceil^{1-\alpha} \right) \geq \frac{1}{\eps},
\]
so that $p_{12}^{\tilde{n}(\eps)} = 1/3$.
An upper bound on $T^0_{1 \to 2}$ is then obtained as (notice that for $\eps$ small enough, $\left\lceil \tilde{b}(\eps) \right\rceil -2 \geq \tilde{n}(\eps)$):
\begin{align}
  \PP\left(T^0_{1 \to 2} \geq \tilde{b}(\eps)\right)
  &=\prod_{k=0}^{\left\lfloor \tilde{b}(\eps)\right\rfloor-2} p_{11}^k 
  \leq \prod_{k=\tilde{n}(\eps)}^{\left\lfloor \tilde{b}(\eps)\right\rfloor-2} p_{11}^k
  \leq \left(\frac{2}{3}\right)^{\left\lfloor \tilde{b}(\eps)\right\rfloor-2-\tilde{n}(\eps)} 
  \nonumber \\
  &\leq \frac{81}{16} \exp \left[\left(\ln \frac32\right) \left(\tilde{C}-C_{\tilde b}\right) |\ln \eps|^{1/(1-\alpha)}\right].\nonumber
\end{align}
The right-hand side goes to zero when $\eps$ goes to $0$, which yields the result for the
asymptotic upper bound $\tilde{b}(\eps)$. This
ends the proof of Lemma~\ref{lem:T120} in the case $\alpha \in (0,1)$.
 
In the case $\alpha=1$, 
for the lower bound, we choose $a(\eps)=f(\eps) \,\eps^{-1/(1+\gamma_\star)}$
for any function $f$ such that $\dps \lim_{\eps \to 0} f(\eps)= 0$.
We have, using \eqref{eq:prod_cas1} for the fourth inequality and denoting by $C$ a positive constant which may change from line to line
\begin{align*}
\PP\Big(T_{1\to 2}^0 \leq a(\eps)\Big) 
& = 1 - \PP\Big(T_{1 \to 2} > \lfloor a(\eps)\rfloor\Big) 
= 1 - \prod_{k=0}^{\lfloor a(\eps)\rfloor-1} p_{11}^k \\
& \leq 1 - \left(p_{11}^{\lfloor a(\eps)\rfloor}\right)^{\lfloor a(\eps)\rfloor}
= 1 - \exp \Big(\lfloor a(\eps)\rfloor \ln \left(p_{11}^{\lfloor
  a(\eps)\rfloor}\right) \Big) 
\leq - \lfloor a(\eps)\rfloor \ln \left(p_{11}^{\lfloor a(\eps)\rfloor}\right) \nonumber \\ 
&= - \lfloor a(\eps)\rfloor \ln \left(1- \frac13 \left( \eps\Xi_{\lfloor a(\eps)\rfloor}
\wedge 1 \right)\right) 
\leq -  \lfloor a(\eps)\rfloor \ln \left(1- \frac13 \eps\Xi_{\lfloor a(\eps)\rfloor}
  \right) \nonumber \\
& \leq  -  \lfloor a(\eps)\rfloor \ln \left(1- C  \eps a(\eps)^{\gamma_\star}
  \right)\nonumber 
\leq  -  \lfloor a(\eps)\rfloor \ln \left(1- C  \eps^{1/(1+\gamma_\star)} f(\eps)^{\gamma_\star}
  \right)\nonumber \\
&\leq C \lfloor a(\eps)\rfloor   \eps^{1/(1+\gamma_\star)} f(\eps)^{\gamma_\star} 
\leq  C f(\eps)^{1+\gamma_\star},\nonumber
\end{align*}
which converges to $0$ as $\eps$ goes to $0$.

We now consider the upper bound. We set 
$b(\eps)=g(\eps) \eps^{-1/(1+\gamma_\star)}$ with $\lim_{\eps \to 0} g(\eps)=\infty$.
In the following, we assume that $g$ grows sufficiently slowly so
that  $\lim_{\eps \to 0} \eps (b(\eps))^{\gamma_\star}=0$. This is not
a restrictive  assumption since the probability $\PP(T^0_{1 \to 2} \geq b(\eps))$ is even lower when the function $g$ goes faster to infinity. Moreover, upon replacing $g(\eps)$ by
$\eps^{1/(1+\gamma_\star)} \lfloor\eps^{-1/(1+\gamma_\star)}g(\eps)\rfloor$, we may assume that $b:(0,1)\to{\mathbb N}$. One has

  \begin{align}
    \PP\left(T^0_{1 \to 2} \geq b(\eps)\right)
    &=\prod_{k=0}^{b(\eps)-2} p_{11}^k 
    = \prod_{k=0}^{b(\eps)-2} \left( 1- \frac13 \left( \eps\Xi_k
    \wedge 1 \right) \right).\nonumber
\end{align}
For $k \leq b(\eps)$, it holds 
$\eps\Xi_k \leq \eps \Xi_{\lfloor b(\eps) \rfloor}$ with the right-hand-side smaller than $C \eps b(\eps)^{\gamma_{\star}}$ by~\eqref{eq:prod_cas1}. This upper-bound goes to zero as $\eps$ goes to zero by assumption. Thus, for $\eps$
sufficiently small,
  \begin{align}
  \PP\Big(T^0_{1 \to 2} \geq b(\eps)\Big)
 & = \prod_{k=0}^{b(\eps)-2} \left( 1-\frac13 \eps\Xi_k
    \right)
    \leq \prod_{k=0}^{b(\eps)-2} \left( 1- C \eps k^{\gamma_\star}
    \right), \label{eq:T12large_cas1_}
  \end{align}
where $C$ is a constant independent of $\eps$.  Then, using the fact that $\eps b(\eps)^{\gamma_\star}$ is smaller than $1/C$ for $\eps$ sufficiently small, we have in this limit
\begin{align*}
  \ln \left( \prod_{k=0}^{b(\eps)-2} \left( 1- C \eps k^{\gamma_\star}
  \right) \right)
  &=
  \sum_{k=0}^{b(\eps)-2} \ln \left(1- C \eps k^{\gamma_\star} \right)
  \leq - C \eps \sum_{k=1}^{b(\eps)-2} k^{\gamma_\star} \\
  &\leq - C \eps \sum_{k=1}^{b(\eps)-2} \int_{k-1}^k x^{\gamma_\star} \, dx
  = - C \eps \int_0^{b(\eps)-2} x^{\gamma_\star} \, dx \\
  & = - \frac{C}{{\gamma_\star}+1}\eps (b(\eps)-2)^{{\gamma_\star}+1} 
  \leq - \frac{C}{{\gamma_\star}+1}\left(1-\frac{2}{b(\varepsilon)}\right)^{{\gamma_\star}+1}g(\eps)^{{\gamma_\star}+1}.
\end{align*}
Using this estimate in~\eqref{eq:T12large_cas1_} leads to the existence of a modified positive constant $C$ such that for $\varepsilon$ small enough,
$\PP(T^0_{1 \to 2} \geq b(\eps)) \leq \exp(-C g(\eps)^{{\gamma_\star}+1})$,
the right-hand side going to~0 as $\eps \to 0$.
This therefore concludes the proof of Lemma~\ref{lem:T120} in the case $\alpha=1$.
\end{adem} 
\begin{rem}
Considering the Equation~\eqref{eq:C2}, one could think of replacing
the assumption~\eqref{eq:beta_1} on $\beta$ by the seemingly weaker
one:
$$\lim_{\eps \to 0} \left(
|\ln \eps| - \beta(\eps) \right)^{\alpha/(1-\alpha)} \exp \left(- \beta(\eps)
\right) = 0.$$
But both conditions are equivalent since
\begin{align*}
|\ln \eps|^{\alpha/(1-\alpha)} \exp \left(- \beta(\eps)
\right) &\leq 1_{\{\beta(\eps) \leq |\ln \eps|/2\}}
2^{\alpha/(1-\alpha)} \left(
|\ln \eps| - \beta(\eps) \right)^{\alpha/(1-\alpha)} \exp \left(- \beta(\eps)
\right)\\
&\quad + 1_{\{\beta(\eps) > |\ln \eps|/2\}} \left(|\ln
  \eps|\right)^{\alpha/(1-\alpha)} \exp \left(- |\ln \eps|/2
\right).
\end{align*}
\end{rem}

In order to prove Lemma~\ref{lem:t3}, we need, as explained in the proof of Proposition \ref{lem:T13}, to ensure that $\tilde{\theta}_n(2)$ remains small when $T^0_{1\to 2}>a(\varepsilon)$ and $N_2\leq \Delta (\varepsilon)$.
\begin{lemma}\label{lem:Delta}
Let us assume that $\alpha \in (0,1]$. Let us consider a non-negative
constant $\Delta$ and a constant $a \geq 1$. Let $\nu_2(n)$
denote the number of visits of state $2$ up to time $n$
included. Then, on the event $\{T^0_{1 \to 2} >
{a}\}$, for any $n$ such that $\nu_2(n)\leq \Delta$, 
\begin{equation}
   \label{eq:majoth2}\tilde \theta_n(2)\leq\begin{cases}\exp\left(\frac{\gamma_\star}{1-\alpha} \lfloor a \rfloor^{1-\alpha} \left( \left(
   \frac{\lceil a+ \Delta \rceil + 1}{\lfloor a \rfloor}
   \right)^{1-\alpha} -  1\right)\right)\mbox{ if }\alpha\in (0,1)\\\left(\frac{\lceil a+\Delta\rceil + 1}{\lfloor a \rfloor}\right)^{\gamma_\star}\mbox{ if }\alpha=1
   \end{cases}.
\end{equation}
\end{lemma}

\begin{adem}[~of Lemma \ref{lem:Delta}]
On the event $\{T^0_{1 \to 2} > {a}\}$, for $n$ such that $\nu_2(n)\leq \Delta$, it holds
\[
\tilde \theta_n(2) \le
\prod_{k=\lfloor a \rfloor +1}^{\lceil a + \Delta \rceil + 1} \left(1 +
\frac{\gamma_\star}{k^\alpha}\right).
\]
Now,
\begin{align*}
\ln \left( \prod_{k=\lfloor a \rfloor +1}^{\lceil a + \Delta \rceil + 1} \left(1 +
  \frac{\gamma_\star}{k^\alpha}\right) \right) 
  & =  \sum_{k=\lfloor a \rfloor +1}^{\lceil a + \Delta \rceil + 1} \ln \left(1 +
  \frac{\gamma_\star}{k^\alpha} \right)
  \leq \sum_{k=\lfloor a \rfloor +1}^{\lceil a + \Delta \rceil + 1}
  \frac{\gamma_\star}{k^\alpha} \\ 
  & \leq \sum_{k=\lfloor a \rfloor +1}^{\lceil a + \Delta \rceil + 1}
 \int_{k-1}^{k}  \frac{\gamma_\star}{x^\alpha} \, dx 
 = \int_{\lfloor a \rfloor}^{\lceil a + \Delta \rceil + 1}
\frac{\gamma_\star}{x^\alpha} \, dx.
\end{align*}
When $\alpha=1$,  the right-hand side is equal to
$\gamma_\star\ln\left(\frac{\lceil a+\Delta\rceil + 1}{\lfloor a
    \rfloor}\right)$, which gives the claimed result. When $\alpha \in
(0,1)$, the  right-hand side is equal to
\begin{align*}
\frac{\gamma_\star}{1-\alpha} \lfloor a \rfloor^{1-\alpha} \left( \left(
   \frac{\lceil a + \Delta \rceil + 1}{\lfloor a \rfloor}
   \right)^{1-\alpha} -  1\right)
\end{align*}
which concludes the proof.
\end{adem}

\begin{adem}[~of Lemma \ref{lem:t3}]
Let $\dps c(\eps)=\left\lceil\frac{ b(\eps) - \tilde b(\eps)}{\Delta(\eps)}\right\rceil-1$.
Using the fact that $N_{2 \to 1} \leq N_2$, and recalling that $({\mathcal F}_{n})_{n \geq 0}$ denotes the filtration generated by $((X_n,\theta_n))_{n \geq 0}$, it holds
\begin{align}
\PP&\left(T^0_{1 \to 2} \in \left(a(\eps),\tilde b(\eps)\right), \, N_2 \leq \Delta(\eps),
\, \sum_{i=1}^{N_{2\to1}}T^i_{1\to2} \geq b(\eps) - \Delta(\eps) - \tilde b(\eps) \right) \nonumber \\
&\leq \PP\left(T^0_{1 \to 2} \in \left(a(\eps),\tilde b(\eps)\right), \, N_2 \leq \Delta(\eps), \, \exists i \in \{1,\ldots,N_{2\to1}\}, \, T^i_{1\to2} \geq \frac{ b(\eps) - \tilde b(\eps)}{\Delta(\eps)}-1 \right) \nonumber \\
&\leq \sum_{l=1}^{\Delta(\eps)} \PP\left(T^0_{1 \to 2} \in \left(a(\eps),\tilde b(\eps)\right), \, N_2 \leq \Delta(\eps), \, N_{2\to1}=l, \, \exists i \in \{1,\ldots,l\}, \, T^i_{1\to2} \geq c(\eps)\right)\nonumber \\
& \leq \sum_{l=1}^{\Delta(\eps)} \sum_{i=1}^{l} \PP\left(T^0_{1 \to 2} \in \left(a(\eps),\tilde b(\eps)\right), \, N_2 \leq \Delta(\eps), \, N_{2\to1}=l, \, T^i_{1\to2} \geq c(\eps)\right) \nonumber \\
& = \sum_{i=1}^{\Delta(\eps)} \PP\left(T^0_{1 \to 2} \in \left(a(\eps),\tilde b(\eps)\right), \, N_2 \leq \Delta(\eps), \, N_{2\to1}\geq i, \, T^i_{1\to2} \geq c(\eps)\right) \nonumber \\
&\leq  \sum_{i=1}^{\Delta(\eps)} \EE\left(1_{\{T^0_{1 \to 2}>a(\eps), N_{2\to1}\geq i\}} \PP\left(N_2\leq \Delta(\eps),T^i_{1\to2} \geq c(\eps) \, \left| \, {\mathcal F}_{\tau^i_{2 \to 1}} \right. \right)\right)\label{eq:majTi12}
\end{align}
where $\tau^i_{2 \to 1}$ is defined
by~\eqref{eq:taui21}. We recall that $\nu_2(n)$ denotes the number of visits of state $2$ up to time $n$ included. On $N_{2\to1}\geq i$, $N_2\geq \nu_2\left(\tau^i_{2 \to 1}+T^i_{1\to2}\right)$ and therefore,  by using the strong Markov property of the chain $((X_n,\t_n))_{n\geq 0}$, we obtain that, on the event $\{N_{2\to 1}\geq i\}$,
\begin{align}
   &\PP\left(N_2\leq \Delta(\eps),T^i_{1\to2} \geq c(\eps)
 \, \left| \, {\mathcal F}_{\tau^i_{2 \to 1}} \right.\right)\notag\\
&\leq \EE\left(\left. 1_{\left\{\nu_2(\tau^i_{2 \to 1}+c(\eps)-2)\leq \Delta(\eps),X_{\tau^i_{2 \to 1}}=1,\ldots, X_{\tau^i_{2 \to 1}+c(\eps)-2}=1\right\}}\PP\left(X_{\tau^i_{2 \to 1}+c(\eps)-1}=1 \, \left| \, {\mathcal F}_{\tau^i_{2 \to 1}+c(\eps)-2} \right. \right)\right| \, {\mathcal F}_{\tau^i_{2 \to 1}}\right)\notag\\
&=\EE\left(\left. 1_{\left\{\nu_2(\tau^i_{2 \to 1}+c(\eps)-2)\leq \Delta(\eps),X_{\tau^i_{2 \to 1}}=1,\ldots, X_{\tau^i_{2 \to 1}+c(\eps)-2}=1\right\}}\Big(1-P_{\theta_{\tau^i_{2 \to 1}+c(\eps)-2}}(1,2)\Big)\,\right| \, {\mathcal F}_{\tau^i_{2 \to 1}}\right).\label{probacond}
\end{align}
We recall that
\[
P_{\theta_n}(1,2)=\frac13 \left(\eps
  \frac{\tilde \theta_n(1)}{\tilde \theta_n(2)} \wedge 1 \right).
\]
On the event $\{T^0_{1 \to 2}>a(\eps)\}$, we have, for $n\geq a(\eps)$,
$\tilde \theta_n(1) \geq \Xi_{a(\eps)}$, so that $P_{\theta_n}(1,2)\geq \frac{\eps  \Xi_{a(\eps)}}{3\tilde \theta_n(2)} \wedge \frac13$.
Since $\Delta(\varepsilon)={\rm O}(a(\varepsilon)^\alpha)$, by Lemma \ref{lem:Delta}, there exist constants $M\in(0,+\infty)$ and $\bar\varepsilon\in(0,1)$ such that for 
\begin{equation}
   \forall\varepsilon\in(0,\bar\varepsilon),\mbox{ on the event }\{T^0_{1 \to 2} > {a}(\eps)\},\;\forall n\mbox{ s.t. }\nu_2(n)\leq \Delta(\eps),\;\tilde{\theta}_n(2)\leq M\label{eq:majotht2}.
\end{equation}
As a consequence, for $\eps\in (0,\bar{\varepsilon})$, on the event $\{T^0_{1 \to 2}>a(\eps)\}\cap\{N_{2\to 1}\geq i\}\cap\{\nu_2(\tau^i_{2 \to 1}+c(\eps)-2)\leq \Delta(\eps)\}$, $\tilde \theta_{\tau^i_{2 \to 1}+c(\eps)-2}(2)\leq M$ and therefore $P_{\theta_{\tau^i_{2 \to 1}+c(\eps)-2}}(1,2)\geq \frac{\eps  \Xi_{a(\eps)}}{3M} \wedge \frac13$.
Since, from~\eqref{eq:XiUB}, 
$\eps  \Xi_{a(\eps)} \leq \exp(-\beta(\eps))$ which goes to zero as $\eps$ goes to~0, we deduce that, up to diminishing $\bar{\varepsilon}$, for any $\varepsilon\in (0,\bar{\varepsilon})$,
\[
P_{\theta_{\tau^i_{2 \to 1}+c(\eps)-2}}(1,2) \geq \frac{\eps  \Xi_{a(\eps)}}{3\tilde M}\mbox{ on the event }\{T^0_{1 \to 2}>a(\eps)\}\cap\{N_{2\to 1}\geq i\}\cap\{\nu_2(\tau^i_{2 \to 1}+c(\eps)-2)\leq \Delta(\eps)\}.
\]With \eqref{probacond}, we deduce that on $\{T^0_{1 \to 2}>a(\eps)\}\cap\{N_{2\to 1}\geq i\}$, 
\begin{align*}
  &\PP\left(N_2\leq \Delta(\eps), \, T^i_{1\to2} \geq c(\eps)
  \, \left| \, {\mathcal F}_{\tau^i_{2 \to 1}}\right.\right)\notag\\
&\leq \EE\left(\left.1_{\left\{\nu_2(\tau^i_{2 \to 1}+c(\eps)-3)\leq
    \Delta(\eps),X_{\tau^i_{2 \to 1}}=1,\ldots, X_{\tau^i_{2 \to
        1}+c(\eps)-3}=1\right\}} \left(1-\frac{\eps  \Xi_{a(\eps)}}{3M}\right)
  \left(1-P_{\theta_{\tau^i_{2 \to 1}+c(\eps)-3}}(1,2)\right) \,\right| \, {\mathcal F}_{\tau^i_{2 \to 1}}\right).
\end{align*}
Iterating the reasoning, we obtain that, on $\{T^0_{1 \to 2}>a(\eps)\}\cap\{N_{2\to 1}\geq i\}$, 
\begin{align*}
\PP&\left(N_2\leq \Delta(\eps), \, T^i_{1\to2} \geq c(\eps)
  \, \left| \, {\mathcal F}_{\tau^i_{2 \to 1}}\right.\right)   
\leq \left(1-\frac{\eps  \Xi_{a(\eps)}}{3M} \right)^{c(\eps)-1}.
\end{align*}
With \eqref{eq:majTi12} and the definition of $c(\varepsilon)$, we deduce that
\begin{align}
& \PP\left(T^0_{1 \to 2} \in \left(a(\eps),\tilde b(\eps)\right), \, N_2 \leq \Delta(\eps),
\, \sum_{i=1}^{N_{2\to1}}T^i_{1\to2} \geq b(\eps) -
\Delta(\eps) - \tilde b(\eps) \right)\notag\\
&\qquad \leq \Delta(\eps)\exp \left(\left(\left\lceil\frac{ b(\eps) - \tilde b(\eps)}{\Delta(\eps)}\right\rceil-2\right) \ln \left( 1 - \frac{\eps\Xi_{a(\eps)}}{3M}
\right)\right)\label{eq:majot3ttal}\end{align}
For $\alpha\in [1/2,1)$, we deduce that
\begin{align*}
   \PP&\left(T^0_{1 \to 2} \in \left(a(\eps),\tilde b(\eps)\right), \, N_2 \leq \Delta(\eps),
\, \sum_{i=1}^{N_{2\to1}}T^i_{1\to2} \geq b(\eps) -
\Delta(\eps) - \tilde b(\eps) \right)\\&\leq \Delta(\eps) \exp \left( - K \frac{1}{\Delta(\eps)} |\ln \eps|^{1/(1-\alpha)}
\eps  \Xi_{a(\eps)}  \right)
\end{align*}
for some positive constant $K > 0$. 
We conclude by \eqref{eq:XiLB} which ensures\begin{align*}
\frac{1}{\Delta(\eps)} |\ln \eps|^{1/(1-\alpha)}
\eps  \Xi_{a(\eps)} 
&\geq C \frac{1}{\Delta(\eps)} |\ln \eps|^{1/(1-\alpha)}
\eps  \exp \left( \frac{\gamma_\star}{1-\alpha}
  a(\eps)^{1-\alpha}\right) \\
&=  C \frac{1}{\Delta(\eps)} |\ln \eps|^{1/(1-\alpha)}
  \exp \left( - \beta(\eps)\right). 
\end{align*}
for $\alpha \in
(1/2,1)$ and
\begin{align*}
\frac{1}{\Delta(\eps)} |\ln \eps|^{1/(1-\alpha)}
\eps  \Xi_{a(\eps)} 
&\geq C \frac{1}{\Delta(\eps)} |\ln \eps|^{2}
\eps  \exp \left( 2 \gamma_\star
  \sqrt{a(\eps)} - \frac{\gamma_\star^2}{2} \ln (a(\eps)) \right) \\
&=  C \frac{1}{\Delta(\eps)} |\ln \eps|^{2}
  \exp \left( - \beta(\eps)\right) \big(|\ln \eps| -
  \beta(\eps)\big)^{-\gamma_\star^2} \\
&\geq  C \frac{1}{\Delta(\eps)} |\ln \eps|^{2- \gamma_\star^2}
  \exp \left( - \beta(\eps)\right).
\end{align*} for $\alpha = 1/2$.

When $\alpha=1$, \eqref{eq:majot3ttal} implies
\begin{align*}
  \PP&\left(T^0_{1 \to 2} \in \left(a(\eps),\tilde b(\eps)\right), \, N_2 \leq \Delta(\eps),
  \, \sum_{i=1}^{N_{2\to1}}T^i_{1\to2} \geq b(\eps) - \Delta(\eps) - \tilde b(\eps) \right)\\
  &\qquad \qquad \leq \Delta(\eps) \exp \left( \left( \frac{ g(\eps) - \tilde g(\eps)}{\Delta(\eps)} \, \eps^{-1/(1+\gamma_\star)} - 2 \right) \ln \left( 1 - \frac{\eps  \Xi_{a(\eps)}}{3M}\right)\right).
\end{align*}
Using the fact that, by \eqref{eq:prod_cas1}, there exists a constant $C$ independent of $\eps$ such that
\[
\eps \Xi_{a(\eps)} \leq C f (\eps)^{\gamma_\star} \eps^{1/(1+\gamma_\star)},
\]
so that the left-hand side goes to zero when $\eps$ goes to zero, we
obtain (the constants $C,C'$ are independent from $\eps$ small enough, and their values may
change from one occurrence to another)
\begin{align*}
\PP&\left(T^0_{1 \to 2} \in \left(a(\eps),\tilde b(\eps)\right), \, N_2 \leq \Delta(\eps),
\, \sum_{i=1}^{N_{2\to1}}T^i_{1\to2} \geq b(\eps) - \Delta(\eps) - \tilde b(\eps) \right)\\
&\qquad \leq C' \Delta(\eps) \exp \left( - C \frac{ g(\eps) - \tilde g(\eps)}{\Delta(\eps)} \, \eps^{-1/(1+\gamma_\star)} \eps  \Xi_{a(\eps)}  \right) \\
&\qquad \leq C'\Delta(\eps)  \exp \left( - C \frac{ g(\eps) - \tilde g(\eps)}{\Delta(\eps)} 
f(\eps)^{\gamma_\star}\right).
\end{align*}\end{adem}

\section{Discussion of the successive exit times of the metastable states}\label{sec:succtimes}

In this section, we consider the scaling of the successive transition
times back and forth between states 1 and 3,  and not only of the
first transition time from 1 to 3. For the sake of conciseness, we do
not provide complete proofs of the results, but only indicate how to
adapt the previous reasoning to the successive exit times.

For the non-adaptive dynamics $\{\overline{X}_n, \, n \ge 0\}$, the
analysis is very easy. Let $\overline{T}_{3\to 1}$ denote the time between $\overline{T}_{1\to 3}$ and the first subsequent return to state $1$ : $\overline{T}_{3\to 1}=\min\{n>\overline{T}_{1\to 3}:\overline{X}_n=1\}-\overline{T}_{1\to 3}$. Of course by symmetry, the asymptotic behavior of $\eps\overline{T}_{3\to 1}$ as $\eps\to 0$ is the same as the one of $\eps\overline{T}_{1\to 3}$ given by Proposition \ref{lem:non-adapt}: $\eps\overline{T}_{3\to 1}$ scales like~$6/\eps$ and converges in distribution to an exponential random variable
with parameter $1/6$. And more generally, all the successive durations needed by the Metropolis-Hastings algorithm to go from one of the extremal states $1$ and $3$ to the other scale like~$6/\eps$.

Let us now discuss the successive exit times of the Wang-Landau
algorithm. We first consider the easier case $\alpha\in[1/2,1)$ in
Section~\ref{sec:1}, which is illustrated by numerical experiments in
Section~\ref{sec:2}. We finally discuss the case $\alpha=1$ in Section~\ref{sec:3}.

\subsection{Successive exit times of the Wang-Landau algorithm for $\alpha\in[1/2,1)$}\label{sec:1}

Setting $n(\eps)=\left(\frac{1-\alpha}{\gamma_\star}\right)^{1/(1-\alpha)}| \ln \eps|^{1/(1-\alpha)}$, one has  $T_{1 \to 3}\sim n(\eps)$ according to Proposition \ref{lem:analytical_adaptive}. Let $T_{3\to 1}$ denote the time between ${T}_{1\to 3}$ and the first subsequent return to state $1$ : ${T}_{3\to 1}=\min\{n>{T}_{1\to 3}:{X}_n=1\}-{T}_{1\to 3}$. To analyse the asymptotic behavior of $T_{3\to 1}$ as $\eps\to 0$, one needs the vector $\tilde{\theta}_{T_{1 \to 3}}$ of unnormalized weights at time $T_{1\to 3}$. One has $\tilde{\theta}_{T_{1 \to 3}}(3)=1+\gamma_\star T_{1\to 3}^{-\alpha}=1+{\rm o}(1)$. By the proof of Proposition \ref{lem:analytical_adaptive} (see in particular Lemma \ref{lem:T120} and \eqref{eq:majotht2}), there is a finite constant $M$ such that $\lim_{\eps\to 0}\PP(\tilde{\theta}_{T_{1 \to 3}}(2)\leq M)=1$. Last, since before time $T_{1\to 3}$, the algorithm stays in state $1$ at least during the time interval $[0,T_{1\to 2}-1]$ and at most during the time interval $[0,T_{1\to 3}-2]$, \begin{equation}\label{eq:contth1}\Xi_{T^0_{1\to 2}-1}\leq \tilde{\theta}_{T^0_{1\to 2}}(1)\leq \tilde{\theta}_{T_{1\to 3}}(1)\leq \Xi_{T_{1\to 3}-2}.\end{equation} For $c\in (1,+\infty)$, choosing $C_a=\left(\frac{1-\alpha}{c\gamma_\star}\right)^{1/(1-\alpha)}$ and $C_b=\left(\frac{c(1-\alpha)}{\gamma_\star}\right)^{1/(1-\alpha)}$, one deduces by Lemma \ref{lem:T120}, Proposition \ref{lem:analytical_adaptive}  and \eqref{eq:equivalent_lnXi_alpha<1},  that 
$$\lim_{\eps\to 0}\PP\left(\frac{1}{c}|\ln \eps|\leq \ln(\tilde{\theta}_{T_{1\to 3}}(1))\leq c|\ln \eps|\right)=1.$$
This means that $\tilde{\theta}_{T_{1\to 3}}(1)$ is approximately of order $\frac{1}{\eps}$. We will perform the analysis of $T_{3\to 1}$, under the simplifying assumption that $\tilde{\theta}_{T_{1\to 3}}(1)\leq \frac{C}{\eps}$ so that, as long as state $1$ has not been reached again after $T_{1\to 3}$, the transition probability from state $2$ to state $1$ remains of order $1$. Then the only difference with the analysis of $T_{1\to 3}$ is that the stepsizes of the Wang-Landau algorithm have been shifted into $\left(\frac{\gamma_\star}{(T_{1\to 3}+n)^\alpha}\right)_{n\geq 1}$. Repeating the analysis performed in the proof of Lemma \ref{lem:T120}, we see that the time $T^0_{3\to 2}$ needed by the algorithm to reach again state $2$ will be of order $n_2(\eps)$ such that
$$\sum_{k=n(\eps)}^{n_2(\eps)}\exp\left(\sum_{j=n(\eps)}^k\frac{\gamma_\star}{j^\alpha}\right)=\mathrm{O}\left(\frac{1}{\eps}\right).$$
This condition gives $n_2(\eps)\sim \left(\frac{2(1-\alpha)}{\gamma_\star}\right)^{1/(1-\alpha)}| \ln \eps|^{1/(1-\alpha)}$. So, repeating the arguments given in the proof of Proposition \ref{lem:analytical_adaptive}, one expects that $T_{3\to 1}$ is of order $\left(\frac{2(1-\alpha)}{\gamma_\star}\right)^{1/(1-\alpha)}| \ln \eps|^{1/(1-\alpha)}$ and that 
$\tilde{\theta}_{T_{1\to 3}+T_{3\to 1}}(2)$ remains bounded. Moreover, one also expects that $\tilde{\theta}_{T_{1\to 3}+T_{3\to 1}}(3)$ is approximately of order $\frac{1}{\eps}$.

At time $T_{1\to 3}+T_{3\to 1}$, one has $\tilde{\theta}_{T_{1\to 3}+T_{3\to 1}}(2)$ bounded uniformly in $\eps$ whereas $\tilde{\theta}_{T_{1\to 3}+T_{3\to 1}}(1)$ and $\tilde{\theta}_{T_{1\to 3}+T_{3\to 1}}(3)$ are both approximately of order $\frac{1}{\eps}$ so that every entry in the transition matrix $P_{{\theta}_{T_{1\to 3}+T_{3\to 1}}}$ but the ones with indices $(1,3)$ and $(3,1)$ are approximately of order $1$. So one expects, that after time $T_{1\to 3}+T_{3\to 1}$ which is of order $\left(1+2^{1/(1-\alpha)}\right)\left(\frac{1-\alpha}{\gamma_\star}\right)^{1/(1-\alpha)}| \ln \eps|^{1/(1-\alpha)}$, the Wang-Landau algorithm has got rid of the initial metastability and moves freely from any of the extremal states $1$ and $3$ to the other one with the only constraint of going through state $2$.

\subsection{Numerical results}\label{sec:2}

\begin{figure}\begin{center}
\includegraphics[width=7.2cm]{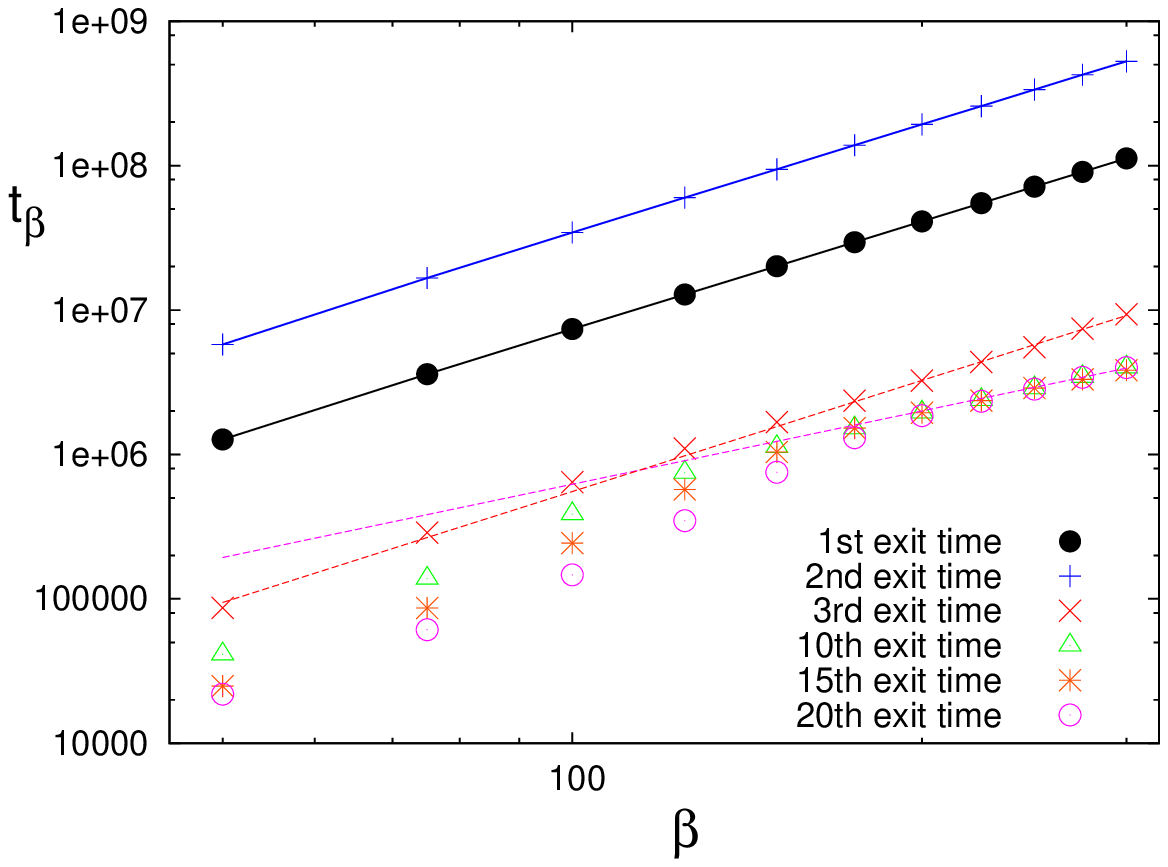}
\includegraphics[width=7.2cm]{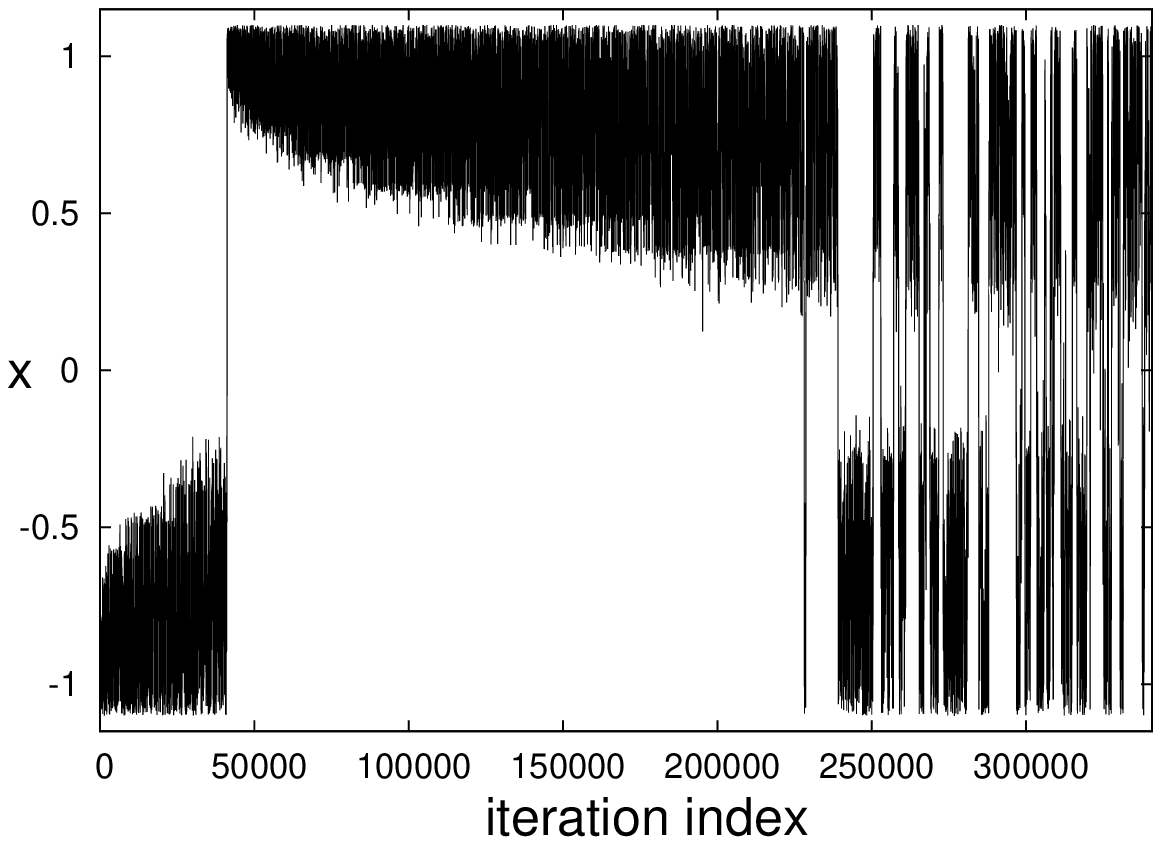}
\caption{\label{fig:many_exits}
  Left: scaling of successive exit times as a function of the inverse temperature (in log-log scale) in the case $\alpha = 0.6$ and $\gamma_\star = 1$. The first three exit times are of the same order of magnitude. All the subsequent exit times have similar orders of magnitudes. The exit times starting from the third one are much smaller than the first two. 
  Right: typical trajectory for $\beta = 15$ when $\alpha = 0.6$ and $\gamma_\star = 1$. Note how the system first explores the two metastability basins before more freely switching from one basin to the other.
}
\end{center}
\end{figure}

The above theoretical results on the scaling of the exit times for a
simple three-state model can be numerically checked for the model
presented in Section~\ref{sec:exit_numerical}. We present in
Figure~\ref{fig:many_exits} the average successive exit times as a
function of the inverse temperature in the case $\alpha = 0.6$ and $\gamma_\star = 1$, as well as a typical trajectory
in order to visualize more clearly the qualitative behavior of the
system. We denote by $t_\beta^k$ the average $k$-th exit time, obtained
by averaging exit times obtained for $M=10^5$ independent
realizations for the smallest values of $\beta$, and a few thousands for the largest values of $\beta$ (the other parameters being the same as in
Section~\ref{sec:numres}, namely $R=1.1$, $d = 22$, $\upsilon = 0.1$). The time $t_\beta^1$ is the first transition time $t_\beta$ introduced in
Section~\ref{sec:exit_numerical}, $t_\beta^2$ is the average
of the first transition time from the value $x_1 = 1$ back to $x_1 = -1$, 
$t_\beta^3$ is the average of the second transition time from the value $x_1=-1$ to $x_1=1$, and so forth. Our numerical results show that 
\[
t_\beta^k \sim C_k \beta^a 
\]
with $a=2.5$ for $k = 1,2,3$, while $a \simeq 1.7$ for
$k \geq 4$. Several conclusions can be
drawn. First, the first two exit times indeed have the same
scaling, as expected from the analysis in the previous section. Moreover, the subsequent exit times (except for the third one) also have the same
scalings, but are much shorter in average than the first two exit
times. They are however still growing with~$\beta$. This
is due to the fact that, in this case which is more complex than the
simple three-state model, some metastability
remains, as illustrated by Figure~\ref{fig:contour} (Right): there are
still  energy (or free energy) barriers to cross, even for the biased potential. 

\subsection{Successive exit times of the Wang-Landau algorithm for $\alpha=1$}\label{sec:3}

In the case $\alpha=1$, by Proposition \ref{lem:analytical_adaptive}, $T_{1\to 3}$ is approximately of order $n(\eps)=\eps^{-1/(1+\gamma_\star)}$. One still has $\tilde{\theta}_{T_{1 \to 3}}(3)=1+\gamma_\star T_{1\to 3}^{-1}=1+{\rm o}(1)$ and $\tilde{\theta}_{T_{1 \to 3}}(2)$ bounded uniformly in $\eps$ small enough. Moreover, \eqref{eq:contth1}, Proposition \ref{lem:analytical_adaptive}, Lemma \ref{lem:T120} and \eqref{eq:prod_cas1} imply that for any function $h$ such that $\lim_{\eps\to 0}h(\eps)=+\infty$, 
$$\lim_{\eps\to 0}\PP\left(h(\eps)^{-\gamma_\star}\eps^{-\gamma_\star/(1+\gamma_\star)}\leq \tilde{\theta}_{T_{1 \to 3}}(1)\leq h(\eps)^{\gamma_\star}\eps^{-\gamma_\star/(1+\gamma_\star)}\right)=1.$$
In particular $\tilde{\theta}_{T_{1 \to 3}}(1)\leq \frac{C}{\eps}$. Now the time $T^0_{3\to 2}$ needed  by the algorithm to reach again state $2$ will be of order $n_2(\eps)$ such that $\sum_{k=n(\eps)}^{n_2(\eps)}\prod_{j=n(\eps)}^{k}\left(1+\frac{\gamma_\star}{j}\right)=\mathrm{O}\left(\frac{1}{\eps}\right)$. With this condition, we deduce that $T^0_{3\to 2}$ and $T_{3\to 1}$ will be approximately of order $\epsilon^{-\frac{1+2\gamma_\star}{(1+\gamma_\star)^2}}$. As a consequence $T_{1\to 3}={\rm o}(T_{3\to 1})$, which we could guess from the explosion as $\alpha\to 1$ of the factor $2^{1/(1-\alpha)}$ appearing in the analysis for $\alpha\in [1/2,1)$. Now, while $\tilde{\theta}_{T_{1\to 3}+T_{3\to 1}}(2)$ remains bounded, $\tilde{\theta}_{T_{1\to 3}+T_{3\to 1}}(1)$ is approximately of order $\eps^{-\gamma_\star/(1+\gamma_\star)}$ while $\tilde{\theta}_{T_{1\to 3}+T_{3\to 1}}(3)$ is approximately of order $\prod_{j=n(\eps)}^{n_2(\eps)}\left(1+\frac{\gamma_\star}{j}\right)$ i.e. of order $\eps^{-\left(\frac{\gamma_\star}{1+\gamma_\star}\right)^2}$. So there remains some metastability preventing the algorithm to move quickly from any of the extremal states to the other one. For instance, the time it will need after $T_{1\to 3}+T_{3\to 1}$ to go back to state $3$ will be approximately of order $\eps^{-(1+2\gamma_\star+2\gamma_\star^2)/(1+\gamma_\star)^3}$ which is intermediate between the orders of $T_{1\to 3}$ and $T_{3\to 1}$. Next, it will take a time of approximate order $\eps^{-(1+2\gamma_\star)/(1+\gamma_\star)^2}$ to go back to state $1$ and the next transition times should be smaller since the orders of $\tilde{\theta}(1)$ and $\tilde{\theta}(3)$ have increased but the shift in the sequence of stepsizes is only multiplied by a constant.

\section*{Acknowledgment}

We thank Eric Moulines for stimulating discussions.  This work is
supported by the French National Research Agency under the grants
ANR-09-BLAN-0216-01 (MEGAS) and ANR-08-BLAN-0218 (BigMC).


\end{document}